\documentclass[11pt]{article}
\usepackage{fullpage,graphicx,psfrag,amsmath,amsfonts,amsthm, amssymb, verbatim}
\usepackage{cite}
\usepackage{cases}
\usepackage{epstopdf}
\usepackage{mdwlist}
\usepackage{bigstrut,array,multirow,tabularx}
\usepackage{xcolor}
\usepackage{subfigure}
\usepackage{dsfont}
\usepackage{bbm,bm}
\usepackage{url}
\usepackage[flushleft]{threeparttable} 
\usepackage{booktabs,caption}
\usepackage{algorithm,algorithmic}
\usepackage[final]{changes}

\usepackage{pifont}

\DeclareMathAlphabet{\mathsf}{OML}{cmbr}{m}{it}

\newtheorem{definition}{Definition}
\newtheorem{theorem}{Theorem}
\newtheorem{lemma}{Lemma}
\newtheorem{corollary}{Corollary}
\newtheorem{proposition}{Proposition}

\newtheorem{remark}{Remark}

\DeclareMathOperator{\tr}{\mathrm{tr}}

\DeclareMathOperator*{\argmin}{argmin}
\DeclareMathOperator*{\argmax}{argmax}

\newcommand{\N}{\mathbbmss{N}}
\newcommand{\R}{\mathbbmss{R}}

\newcommand{\B}[1]{\mathbf{#1}}

\newcommand{\BB}[1]{\ensuremath{\bm{#1}}}

\newcommand{\IM}[1]{\B{I}_{#1}}

\newcommand{\Vnorm}[1]{\left\|{#1}\right\|}

\newcommand{\InP}[2]{\left\langle #1,#2 \right \rangle}

\newcommand{\diag}[1]{\mathrm{diag}\left(#1\right)}

\newcommand{\xv}{\BB{x}}
\newcommand{\vv}{\BB{v}}
\newcommand{\uv}{\BB{u}}

\newcommand{\wv}{\BB{w}}
\newcommand{\wtv}{\tilde{\BB{w}}}
\newcommand{\yv}{\BB{y}}

\newcommand{\nv}{\BB{e}}
\newcommand{\anv}{\BB{\xi}}
\newcommand{\Em}{\BB{\Xi}}
\newcommand{\nvk}{\B{E}}

\newcommand{\cX}{\BB{\phi}}
\newcommand{\cY}{\BB{\varphi}}
\newcommand{\CX}{\BB{\Phi}}
\newcommand{\CY}{\BB{\Psi}}
\newcommand{\augv}{\BB{\vartheta}}
\newcommand{\Ym}{\B{Y}}
\newcommand{\Zm}{\B{Z}}

\newcommand{\Bm}{\B{B}}
\newcommand{\bigB}{\B{\mathcal{B}}}
\newcommand{\Cm}{\B{C}}
\newcommand{\Am}{\B{A}}
\newcommand{\Wm}{\B{W}}
\newcommand{\Vm}{\B{V}}
\newcommand{\Vmp}{\B{V}_{\perp}}
\newcommand{\Um}{\B{U}}
\newcommand{\Eig}{\B{\Lambda}}
\newcommand{\Eeig}{\B{\Lambda}}
\newcommand{\Eeigp}{\Eeig_\perp}
\newcommand{\Wtm}{\tilde{\B{W}}}
\newcommand{\Rm}{\B{R}}
\newcommand{\Cov}[1]{\B{\Sigma}_{#1}}
\newcommand{\angl}[2]{\theta_{#1}\left(#2\right)}
\newcommand{\ires}{\textsf{r}_\mathrm{init}}
\newcommand{\dres}{\textsf{r}_\mathrm{des}}
\newcommand{\Cheb}[2]{T_{#1}\left(#2\right)}
\newcommand{\Pv}{\mathrm{P}_{\vv_1}}
\newcommand{\PPv}{\mathrm{P}^{\perp}_{\vv_1}}

\newcommand{\eig}{\lambda}
\newcommand{\gap}{\Delta}
\newcommand{\reg}{\gamma}

\newcommand{\Mnorm}[2]{\left\Vert{#1}\right\Vert_{#2}}
\newcommand{\norm}[1]{\left|{#1}\right|}

\definechangesauthor[color=red]{VM}

\begin{document}

\title{
		 Noisy Accelerated Power Method for Eigenproblems with Applications
}

\author{
	{Vien V. Mai and  Mikael Johansson}{\thanks{V. V. Mai and M. Johansson are with the Department of Automatic Control, School of Electrical Engineering and Computer Science, Royal Institute of Technology (KTH), SE-100 44 Stockholm, Sweden. Emails: \tt\small\{maivv, mikaelj\}@kth.se.}%
	}       
}

\maketitle

\begin{abstract}	
	This paper \replaced{introduces an efficient}{studies \replaced{a practically and theoretically efficient}{an efficient and provable}} algorithm for finding the \deleted{top}dominant generalized eigenvectors of a pair of symmetric matrices. Combining tools from approximation theory and convex optimization, we \replaced{develop}{present} a simple scalable algorithm with strong theoretical \replaced{performance guarantees}{guarantee}. More precisely, the algorithm \replaced{retains}{admits} the simplicity of the well-known \emph{power method} \replaced{but}{and} enjoys the\deleted{same} asymptotic iteration complexity of the powerful \emph{Lanczos method}. Unlike these classic \replaced{techniques}{methods}, our algorithm is designed to decompose \added{the overall} problem into a series of subproblems that \added{only} need \added{to be}\deleted{only be} solved approximately. \replaced{The combination of good initializations, fast iterative solvers, and}{Good initializations when combined with} appropriate error control in solving the subproblems lead to a \emph{linear} running time in the input sizes compared to \added{the} \emph{superlinear} time \replaced{for}{as in} the traditional methods. The improved running time \replaced{immediately}{readily} offers acceleration for several applications. 
	As an example, we demonstrate how the proposed algorithm can be used to accelerate canonical correlation analysis, which is a fundamental statistical tool for learning of a low-dimensional representation of high-dimensional objects. Numerical experiments on real-world data sets confirm that our approach yields significant improvements over the current state-of-the-art.
\end{abstract}

\section{Introduction}
Computing the leading eigenvectors of a matrix is one of the most fundamental problems in scientific computing. The problem has been studied \replaced{extensively}{extremely well} in a classical setting and is used in countless applications.  
For generalized eigenvector problems, standard iterative methods  such as the \emph{power method} and the Lanczos method require\deleted{performing} multiple matrix-vector products of the form $\Bm^{-1}\xv$ at each iteration \cite[Section~9.2.6]{Sad11}, which becomes prohibitively expensive \replaced{in high dimensions}{for large-scale data sets}. The poor scalability of the standard techniques make them unfit for many  applications in machine learning and data science.  In these applictions, the ability to deal efficiently with large-scale data is the main concern, while machine precision accuracy is no longer essential since the problem data is often uncertain and subject to noise.  Novel algorithms for scalable and approximate eigenvalue computations are therefore urgently needed.

In this paper, we consider the problem of finding the solution \added{$(\wv, \lambda)$} to the following equation 
\begin{align}\label{eq:EP}
\B{A}\wv = \lambda \B{B} \wv
\end{align}
\deleted{with respect to $\wv$ and $\lambda$,}where $\B{A}$ \replaced{and}{,} $\B{B}$ are \added{given} real, symmetric matrices and $\B{B}$ is positive definite.
\replaced{Such eigenproblems appear in a wide range of machine learning tasks, such as}{The problems in this class range from} finding a set of directions in the
data-embedding space \replaced{which contains}{containing} the maximum amount of variance \deleted{in the
	data}(principal components analysis), \deleted{to}finding a hyperplane that separates
two classes of data \deleted{minimizing a certain cost function}(Fisher discriminant),
\replaced{and}{or} finding correlations between two different representations of the same data
(canonical correlation analysis), among many others.

Let $\B{X}, \B{Y} \in \R^{n \times d}$ be two \replaced{(potentially large)}{large} sets of data\deleted{points} with \deleted{the}empirical covariance matrices $\Cov{11}=\frac{1}{n}\B{X}^\top \B{X}$, $\Cov{22}=\frac{1}{n}\B{Y}^\top\B{Y}$, and  $\Cov{12}=\frac{1}{n}\B{X}^\top\B{Y}$. 
The canonical correlation analysis (CCA) problem \replaced{is concerned}{concerns} with finding the features $\B{x}$ and $\B{y}$ that best encapsulate the similarity of $\B{X}$ and $\B{Y}$. In particular, CCA maximizes the empirical \added{cross-}correlation between $\B{X}$ and $\B{Y}$ by solving the following problem
\begin{align}\label{P:CCA}
\begin{aligned}
& \underset{\B{x}, \B{y}}{\text{maximize}}
& & 
\xv^\top\Cov{12}\yv
\\
& \text{subject to}
&& \xv^\top\Cov{11}\xv= \yv^\top\Cov{22}\yv=1.
\end{aligned}
\end{align}
 Generalized eigenvalue \added{(GEV)} problems, on the other hand, \deleted{(GEV)}compute features that maximize discrepancies between the data sets as measured by the the quantities
\begin{align*}
\underset{\B{x}}{\max} \frac{\xv^\top\Cov{11}\xv}{\xv^\top\Cov{22}\xv} \quad \text{and} \quad
\underset{\B{y}}{\max} \frac{\yv^\top\Cov{22}\yv}{\yv^\top\Cov{11}\yv}.
\end{align*}
Both problems can be reduced to performing \replaced{principal}{principle} component analysis (PCA) on complicated matrices, e.g., $\Cov{22}^{-1/2}\Cov{12}^{\top}\Cov{11}^{-1}\Cov{12}\Cov{22}^{-1/2}$ for CCA and $\Cov{22}^{-1/2}\Cov{11}\Cov{22}^{-1/2}$ for GEV. However, this approach is not scalable, as the formation of $\Cov{11}^{-1/2}$ and $\Cov{22}^{-1/2}$ becomes prohibitively expensive for large data sets.

\subsection{Related Work}
While there have been \replaced{significant recent progress in}{several attempts at} understanding PCA, efficient algorithms and theoretical guarantees remain limited for generalized eigenvalue problems and CCA.
Below, we summarize recent work in optimization and approximation theory that have influenced the developments in this paper. 
\subsubsection{Stochastic optimization for PCA}
Over the past few years, the challenges of dealing with huge data sets have inspired 
the development of novel optimization algorithms for empirical risk minimization (see, e.g., \cite{BCN16} and references therein). 
Leveraging ideas from stochastic optimization, several scalable algorithms have also been proposed for different eigenvalue problems\deleted{stochastic methods have been proposed for the eigenvalue world} \cite{GHJ16,Sha15,Sha16,JJK17,SHM17}. These algorithms try to combine the low iteration cost of the stochastic methods and the high accuracy of the standard deterministic techniques. Notably, borrowing the idea of \emph{variance reduction}~\cite{JZ13}, the author in \cite{Sha15} studied a variant of Oja's algorithm \cite{Oja82}, giving the first \emph{linearly convergent} algorithm for stochastic PCA. The authors of \cite{GHJ16} improved the result further by using the well-known \emph{shift-and-invert} technique from numerical linear algebra~\cite{GV12,Sad11}. Inspired by the \emph{momentum} method as known as the \emph{Heavy ball} method in the optimization literature~\cite{Pol87}, an accelerated stochastic PCA was proposed in \cite{SHM17}. 

\subsubsection{Approximation theory and acceleration}
Approximation theory is ubiquitous in machine learning, optimization, and numerical linear algebra. This is because many problems in these areas can be reduced to finding a faster way to compute primitives  such as $x^s$, $x^{-1}$, or $e^x$~\cite{SV14}. For example, for the solution to linear systems of equations {of} the form $Ax=b$ with $A\in \R^{n \times n}$, Chebyshev iteration improves the iteration complexity from $\mathcal{O}\left(n\kappa\log\frac{1}{\epsilon}\right)$ of the Richardson iteration to $\mathcal{O}\left(n\sqrt{\kappa}\log\frac{1}{\epsilon}\right)$, where $\kappa$ is the condition number of $A$ \cite[Chapter~11]{GV12}. This improved rate is similar to the speed-ups by Nesterov's fast gradient method over the classical gradient descent \cite{Nes04}, which was recently discussed in  terms of polynomial approximations in~\cite{SAB16}. For eigenproblems, the Lanczos method can be considered as an accelerated version of the power method and improves the iteration complexity from  $\mathcal{O}\left(\frac{1}{\gap}\log\frac{1}{\epsilon}\right)$ to  $\mathcal{O}\left(\frac{1}{\sqrt{\gap}}\log\frac{1}{\epsilon}\right)$. Here, $\gap$ is the relative difference between the two eigenvalues of largest magnitude of the relevant matrix, which tends to be very small for many practical data sets. Although polynomial approximation is not explicit in the Lancoz method, its convergence proof is heavily based on Chebyshev polynomials~\cite{Sad11}. We remark that the notion of \emph{noisy} Chebyshev iterations was studied  for linear system of equations in \cite{GO88}; however, the proof technique is not directly applicable to eigenproblems. 

\subsubsection{Canonical correlation analysis}
Recently, a series of papers have appeared which establish fundamental theoretical guarantees and propose efficient algorithms for the CCA problem \cite{CKL09,MLF15,GJK16,WWG16}. One of the first scalable algorithm for solving CCA is \texttt{AppGrad}~\cite{MLF15}, which is an alternating least-squares method followed by projection steps. However, the algorithm is only guaranteed to converge locally.  Inspired by the noisy power method \cite{HP14}, the authors in \cite{GJK16} proposed \texttt{CCALin}, a globally convergent algorithm designed to efficiently approximate the power method by breaking the problem into several subproblems, each of them is solved inexactly in an controllable manner using some of the most advanced convex optimization solvers. The main drawback of CCALin is that it has the same iteration complexity $\mathcal{O}\left(\frac{1}{\gap}\log\frac{1}{\epsilon}\right)$ as the power method. Recently, an accelerated rate for CCA has been achieved by the use of carefully tuned iterative methods~\cite{WWG16,AZL17}.\deleted{ based on the techniques developed in \cite{GHJ16}} 
{In particular, for solving top-1 CCA problems, the shift-and-invert (\texttt{SI}) \cite{WWG16} and the \texttt{LazyCCA} \cite{AZL17} methods can achieve a so-called doubly-accelerated rate, i.e., having running time depends on $\sqrt{\tilde{\kappa}}$ and $1/\sqrt{\gap}$ (see, Table~\ref{table:cca}). For top-$k$ CCA problems, only \texttt{LazyCCA}, which relies on the shift-and-invert framework coupled with an extensive analysis, can achieve doubly-accelerated rate. In this work, we aim to achieve such doubly-accelerated rate based on a simple modification of the power method. Note that to establish such improved convergence bounds, \texttt{LazyCCA}, \texttt{SI}, and our method require more information about the CCA problem than \texttt{CCALin}.}

\begin{table}[t]
	\caption{
		Summary of different algorithms solving CCA
	}
	\label{table:cca}
	\begin{center}
		\begin{threeparttable}
			\begin{tabular}{c|c|c|c}
				\hline \hline
				Problem & Algorithm & LS solver & Time complexity 
				\bigstrut\\
				\hline\hline
				\multirow{7}{*}{CCA-1} 
					&\texttt{Appgrad}\cite{MLF15} 
					& GD  
					& $\mathcal{O}\left(\frac{nd\kappa}{\gap}\log\frac{1}{\epsilon}\right)$\tnote{*}
					\\[0.75ex]	
					\cline{2-4}
					&\texttt{CCALin}\cite{GJK16} 
					& SVRG
					&$\mathcal{O}\left(\frac{d\left(n+\sqrt{n\tilde{\kappa}}\right)}{{\gap}}\log\frac{1}{\epsilon}\right)$ 
					\\[0.75ex] 
					\cline{2-4}
					&\texttt{ALS-VR}\cite{WWG16} 
					& SVRG
					& $\mathcal{O}\left(\frac{d\left(n+\tilde{\kappa}\right)}{\gap^2}\log^2\frac{1}{\epsilon}\right)$ 
					\\[0.75ex]
					\cline{2-4}
					&{\texttt{SI}\cite{WWG16}}
					& ASVRG  
					&$\mathcal{O}\left(\frac{dn^{3/4}\sqrt{\tilde{\kappa}}}{\sqrt{\gap\norm{\eig_1}}}\log^2\frac{1}{\epsilon}\right)$ 
					\\[0.75ex]
					\cline{2-4}
					&{\texttt{LazyCCA}}\cite{AZL17}
					& ASVRG  
					&$\mathcal{O}\left(\frac{d\left(n+\sqrt{n\tilde{\kappa}}\right)}{\sqrt{\gap}}\log\frac{1}{\epsilon}\right)$ 
					\\[0.75ex]
					\cline{2-4}
					&{\texttt{\textbf{NAPI}}}
					& ASVRG  
					&$\mathcal{O}\left(\frac{d\left(n+\sqrt{n\tilde{\kappa}}\right)}{\sqrt{\gap}}\log\frac{1}{\epsilon}\right)$ 
				\\
				\hline\hline
				\multirow{6}{*}{CCA-$k$} 
					&\texttt{Appgrad}\cite{MLF15} 
					& GD  
					& $\mathcal{O}\left(\frac{ndk\kappa}{\gap_k}\log\frac{1}{\epsilon}\right)$\tnote{*}
					\\[0.75ex]	
					\cline{2-4}
					&\texttt{CCALin}\cite{GJK16} 
					&SVRG
					&$\mathcal{O}\left(\frac{dk\left(n+\sqrt{n\tilde{\kappa}}\right)}{{\gap_k}}\log\frac{1}{\epsilon}\right)$ 
					\\[0.75ex]
					\cline{2-4}
					&{\texttt{LazyCCA}}\cite{AZL17}
					& ASVRG  
					&$\mathcal{O}\left(\frac{dk\left(n+\sqrt{n\tilde{\kappa}}\right)}{\sqrt{\gap_k}}\log\frac{1}{\epsilon}\right)$ 
					\\[0.75ex]
					\cline{2-4}
					&{\texttt{\textbf{NAPI}}}
					& ASVRG  
					&$\mathcal{O}\left(\frac{dk\left(n+\sqrt{n\tilde{\kappa}}\right)}{\sqrt{\gap_k}}\log\frac{1}{\epsilon}\right)$ 
				\\
				\hline
			\end{tabular}
			\begin{tablenotes}
				\item[*] Local convergence.
				\item[**] $\kappa$ and $\tilde{\kappa}$ are condition numbers, see Section~IV for the formal definitions.
			\end{tablenotes}
		\end{threeparttable}
	\end{center}
\end{table}

\subsection{Contributions}
We propose a \emph{linearly convergent} algorithm for the GEV and CCA problems which combines the simplicity of the power method with the asymptotic performance of the Lanczos method. 
Our method---called \texttt{NAPI}---builds on an acceleration technique inspired by polynomial approximation and an efficient \emph{inexact} matrix-vector product computation via advanced convex optimization {methods}. Our work is motivated by the work in \cite{GO88,GJK16,SHM17}, and the main contributions are summarized as follows:
\begin{itemize}
	\item We explicitly characterize the convergence properties of \texttt{NAPI}. Our result matches the asymptotic iteration complexity of the Lanczos method but with much lower computational cost per iteration. In particular, \texttt{NAPI} attains the worst-case complexity bound $\mathcal{O}\left(\frac{1}{\sqrt{\gap}}\log\frac{1}{\epsilon}\right)$ which is not improvable in terms of the eigenvalue gap $\Delta$. While the Lanczos method can not operate in a stochastic setting, each iteration of \texttt{NAPI} involves solving a least squares problem inexactly, hence can be handled very efficiently by some of the most advanced convex optimization solvers. By appropriately controlling the errors in solving the least squares problem, and by using a well-designed \emph{initialization}, \texttt{NAPI} obtains a \emph{linear} running time in the input-size  which allows it to operate on large-scale data sets.
	\item Our result can be used to obtain improved running time for several downstream applications.  We demonstrate how \texttt{NAPI} can be applied to the CCA problem, resulting in practical and theoretical improvements over the state-of-the-art algorithms (see, Table~\ref{table:cca}). Experiments on real-world data sets confirm the effectiveness of our method.
\end{itemize}

\subsection{Notation}
The symbols $\R$ and $\N$ denote the set of real and nonnegative natural numbers, respectively.  We use $\Vnorm{\cdot}$ to denote the Euclidean norm for vectors and the spectral norm for matrices. For a symmetric positive definite matrix $\Bm$, we denote $\xv^\top\Bm\yv$ as the $\Bm$-inner product of $\xv$ and $\yv$ and $\Mnorm{\xv}{\Bm}=\sqrt{\xv^\top\Bm\xv}$ as the $\Bm$-norm of $\xv$. {For a matrix $\B{X}$, we denote by $\Vnorm{\B{X}}_{\Bm}=\Vnorm{\Bm^{1/2}\B{X}}$ and $\Vnorm{\B{X}}_{\Bm,\rm{F}}=\Vnorm{\Bm^{1/2}\B{X}}_{\rm{F}}$ the spectral and Frobenius norms, respectively.} The number of nonzero elements of a matrix $\Am$ is denoted by $\mathrm{nnz}\left(\Am\right)$.

\subsection{Preliminaries}
Principal angles are a useful tool for studying the convergence of subspaces in iterative eigenvalue methods. We rely on the following definition:
\begin{definition}[Principal Angles \cite{ZK13}]
Let $\mathcal{X}$ and $\mathcal{Y}$ be subspaces of $\R^d$ of dimension at least $k$. The principal angles $0 \leq \theta^{(1)} \leq \theta^{(2)} \cdots \leq \theta^{(k)} \leq \pi/2$ between $\mathcal{X} $ and $\mathcal{Y}$  are defined recursively via:
\begin{align*}
	\cos \theta^{(i)} \left(\mathcal{X},\mathcal{Y}\right)
	=
		\max_{\xv \in \mathcal{X}}
		\max_{\yv \in \mathcal{Y}}\; 		
		\xv^\top \Bm \yv
	=	
		\norm{\xv_i^\top \Bm \yv_i},
\end{align*}
\replaced{subject to}{where} $\Vnorm{\xv}_{\Bm}=\Vnorm{\yv}_{\Bm}=1$, $\xv^\top\Bm \xv_j=0$, $\yv^\top\Bm \yv_j=0$, $j=1,\ldots, i-1$. The vectors $\{\xv_1,\xv_2, \ldots, \xv_k\}$ and $\{\yv_1,\yv_2, \ldots, \yv_k\}$ are called the principal vectors. 
For matrices $\B{X}$ and $\B{Y}$, we use $\theta^{(i)}\left(\B{X},\B{Y}\right)$ to denote the $i$th principal angle between their ranges.
\end{definition}

Since we \added{are} only \replaced{interested}{interest} in the largest principal angle, we \replaced{make \deleted{a}slight abuse of notation and write}{slightly abuse the notation by writing} $\angl{}{\B{X},\B{Y}}$ instead of $\theta^{(k)}\left(\B{X},\B{Y}\right)$.
\begin{proposition}\label{prop:sin:cos:tan}
Let $\B{X}$ and $\B{Y}$ be orthogonal bases for subspaces $\mathcal{X}$ and $\mathcal{Y}$ w.r.t  $\B{B}$, respectively. Let $\B{X}_{\bot}$ be an orthogonal basis w.r.t $\B{B}$ for the orthogonal complement of $\mathcal{X}$. Then, 
\begin{align*}
	\cos\angl{}{\B{X},\B{Y}}= \sigma_{k}\left(\B{X}^\top\B{B}\B{Y}\right) \,\, \text{and} \,\,\,
	\sin\angl{}{\B{X},\B{Y}} = \Vnorm{\B{X}_{\bot}^{\top}\B{B}\B{Y}}.
\end{align*}
If, in addition, $\B{X}^\top\B{B}\B{Y}$ is invertible, then
\begin{align*}
	\tan\angl{}{\B{X},\B{Y}} = \Vnorm{\B{X}_{\bot}^{\top}\B{B}\B{Y} \left(\B{X}^\top\B{B}\B{Y}\right)^{-1}}.
\end{align*}
\end{proposition}

The following lemma presents a nice connection between the principal angles and the distances between subspaces\deleted{and \deleted{it} will be used \replaced{in the proof of our main result.}{later to derive the main results of our work.}}. 
The proof of this lemma when $\Bm=\IM{}$ can be found in {\cite[Section~II.4]{SS90}}; its generalization to general positive definite matrices is straightforward, hence omitted here.
\begin{lemma}[Distance Between Subspaces]\label{lem:angle:space}
Let $\B{X}$ and $\B{Y}$ be orthogonal bases for subspaces $\mathcal{X}$ and $\mathcal{Y}$ w.r.t $\B{B}$, respectively. Let $\B{Q}$ be an orthogonal matrix, then 
\begin{align*}
		\Vnorm{\B{X}\B{X}^\top\Bm-\B{Y}\B{Y}^\top\Bm} &= \sin\angl{}{\B{X},\B{Y}}, 
		\\		
		\min_{ \B{Q}^\top\B{Q}=\IM{k}}\Vnorm{\Bm^{1/2}\B{X}-\Bm^{1/2}\B{Y}\B{Q}} &\leq 2 \sin\frac{\angl{}{\B{X},\B{Y}}}{2}.
\end{align*}
\end{lemma}

\section{Polynomial Approximation, Least Squares Solvers and Eigenproblems}

\subsection{What does polynomial approximation have to do with eigenvalues?}
The \emph{power method}~\cite{GV12} is one of the oldest techniques for computing the largest eigenvalue of a matrix, and the conceptual starting point for many contemporary numerical algorithms. It works as follows: choose an initial vector $\xv_0 \in \R^d$ uniformly at random from the unit sphere and multiply it repeatedly by the symmetric matrix $\Am \in \R^{d\times d}$. This generates a sequence of iterates $\{\xv_s\}$, related via $\xv_s = \Am\xv_{s-1}=\Am^s\xv_0$.
Under certain conditions, this sequence, normalized appropriately, will converge to the leading eigenvector of $A$, i.e, the one associated with the eigenvalue of largest magnitude. 

 \replaced{The power method can be viewed as applying }{It can be seen that the power method is equivalent to an application} of a \deleted{very basic}degree-$s$ polynomial, namely $p_s\left(z\right)=z^s$, to $\Am$, and then multiplying it by the initial vector $\xv_0$. The polynomial $p_s\left(z\right)$ evaluates $1$ to $1$ and moves every $z$ with $\norm{z} <1$  \added{closer} to $0$ at the rate determined by $\norm{z}$. Let $\Am= \Um\Eig \Um^\top$ be the eigenvalue value decomposition of $\Am$, where $\Eig = \diag{{\eig_1},{\ldots},{\eig_d}}$ is a matrix of eigenvalues satisfying $\norm{\eig_1}> \norm{\eig_2} \geq\ldots \geq\norm{\eig_d}$\replaced{. Then, }{, then} the polynomial  $p_s\left(z\right)$ applied to $\Am$ is simply
\begin{align*}
 	p_s\left(\Am\right)
 	&= 
 		\Um
 		\begin{bmatrix}
			p_s\left(\eig_1\right) & 0 & \ldots & 0 \\[0.3em]
			0 & p_s\left(\eig_2\right) & \ldots & 0  \\[0.3em]
			\vdots & \ddots & \ddots  & 0				\\[0.3em]
			0 & \ldots & 0 & p_s\left(\eig_d\right)  		
		\end{bmatrix}
 		\Um^\top
 	= 
 		\Um 
 		\begin{bmatrix}
			 \eig_1^s & 0 & \ldots & 0 \\[0.3em]
			0 & \eig_2^s & \ldots & 0  \\[0.3em]
			\vdots & \ddots & \ddots  & 0				\\[0.3em]
			0 & \ldots & 0 & \eig_2^s  		
		\end{bmatrix}
 		\Um^\top.
\end{align*}
Since $\norm{\eig_1} > \norm{\eig_2}$, {$\Am^s/\eig_1^s$ converges to $\uv_1\uv_1^\top$}, where $\uv_1$ is the first column of $\Um$. This means that an estimate \replaced{of}{to} $\uv_1$ can be found by simply picking the first column of $\Am^s$ and normalizing it to have a unit norm. \replaced{If we suppose that}{Suppose that} $\norm{\eig_2}=\left(1-\gap\right)\norm{\eig_1}$, then \replaced{it will take roughly}{power method takes about} $s=\mathcal{O}\left( \frac{1}{\gap} \log\frac{1}{\epsilon}\right)$ iterations \replaced{until $\vert \lambda_2\vert ^s\leq \epsilon \vert \lambda_1\vert^s$}{to approximately find $\uv_1$ (and $\eig_1$)}. \replaced{Since }{However, the main drawback of the power method is that} the eigengap $\gap$ tends to be very small for practical large-scale matrices, \replaced{the power method converges very slowly}{which makes it converges slowly}. \replaced{Fortunately, the next proposition indicates that it is possible to }{We wish} to improve the running time from $\mathcal{O}\left( \frac{1}{\gap} \log\frac{1}{\epsilon}\right)$ to $\mathcal{O}\left( \frac{1}{\sqrt{\gap}} \log\frac{1}{\epsilon}\right)$ by a simple modification of the power method. 
\deleted{Fortunately, this improved running time is indeed possible thanks to the following powerful result.}
\begin{proposition}{\cite[Theorem~3.3]{SV14}}\label{prop:poly:existence}
For any positive integers $s$ and $d$, there is a degree-$d$ polynomial $p_{s,d}\left(x\right)$ satisfying 
\begin{align*}
	\sup_{x \in \left[-1,1\right]} \norm{p_{s,d}\left(x\right) - x^s} \leq 2e^{-d^2/2s}.
\end{align*}
\end{proposition}
The result implies that \added{it is possible to approximate $x^s$ to} any \added{accuracy} $\delta >0$ \added{i.e.}
$$\sup_{x \in \left[-1,1\right]} \norm{p_{s,d}\left(x\right) - x^s} \leq \delta.$$ 
\added{using a polynomial $p_{s,d}(x)$ of degree} $d =  \lceil \sqrt{2\ln\left(2/\delta\right)s} \rceil $. 
\added{This suggests that  it should be possible to develop an approximate power iteration which reduces the number of matrix multiplications by $A$ from $s$ to $\sqrt{s}$.}  
\deleted{Thus, Proposition~\ref{prop:poly:existence} allows to obtain a polynomial approximating of $x^s$ with degree roughly $\sqrt{s}$.}Of course, the existence result in Proposition~\ref{prop:poly:existence} is only useful if there is an efficient way to construct the approximating polynomial. The proof of the proposition relies on Chebyshev polynomials, which are ubiquitous in numerical optimization. 
We also note that the result above is essentially optimal in the sense that  polynomial approximations of $x^s$ over $\left[-1,1\right]$ require degree $\Omega(\sqrt{s})$~\cite[Chapter~5]{SV14}.

\subsection{Stochastic optimization of least squares}

Over the past few years, there has been a great attention in developing stochastic algorithms for solving \emph{empirical risk minimization} problems in machine learning. Specifically, one is interested in the unconstrained minimization of
\begin{align}\label{eq:fs}
	f\left(x\right):=\frac{1}{n}\sum_{i=1}^{n} f_i\left(x\right),
\end{align}
where $f_1,\ldots,f_n$ is $L$-smooth and convex functions, and $f$ is $\mu$-strongly convex function. Solving this problem using standard first order methods such as gradient descent requires $\mathcal{O}\left(n\kappa\log\frac{1}{\epsilon}\right)$ passes over the data set to achieve an $\epsilon$-optimal solution, where \added{the condition number} $\kappa = L/\mu$ \deleted{is the \emph{condition number} which}can be as large as $\Omega\left(n\right)$ in machine learning applications \cite{Bub15}. This running time can be further improved to $\mathcal{O}\left(n\sqrt{\kappa}\log\frac{1}{\epsilon}\right)$ if one uses Nesterov's accelerated scheme \cite{Nes04}. Many \emph{linearly convergent} \added{stochastic} methods\added{, such as SDCA \cite{SZ13a}, SAGA \cite{DBL14}, or SVRG \cite{JZ13},} have been introduced and \replaced{shown to outperform the}{outperformed} standard first order methods under mild assumptions. These  methods only access to the gradient of one individual component function $f_i$ at each step of the algorithm, instead of the full gradient as in standard first order methods. This results in an improved computational complexity \replaced{of}{, which requires only} $\mathcal{O}\left(\left(n+\kappa\right)\log \frac{1}{\epsilon}\right)$ passes over the data set to achieve an $\epsilon$-optimal solution in expectation. \replaced{When these methods are combined with Nesterov acceleration, the expected complexity becomes}{Note that one can obtain a Nesterov's accelerated version of these algorithms that would need only} $O\left((n+\sqrt{n\kappa})\log(1/\epsilon)\right)$ (see, e.g., \cite{FGKS15b} and \cite{All16}).

{In this work, we approximately compute} the matrix-vector product $\Bm^{-1}\BB{b}$ \replaced{by minimizing the least-squares loss}{via solving the following least squares problem} 
\begin{align*}
	f\left(\wv\right) := \frac{1}{2}\wv^\top \B{B}\wv - \wv^\top\BB{b},
\end{align*} 
over $\wv$. This problem has the unique solution is $\wv^\star =\Bm^{-1}\BB{b}$. In our applications of interest, $\Bm$ is usually an empirical covariance matrix of the form $\frac{1}{n}\sum_{i=1}^n \xv_i\xv_i^\top$, hence $f\left(\wv\right)$ can be expressed as a finite-sum as in~(\ref{eq:fs}) with 
\begin{align*}
f_i(\wv) &:=  \frac{1}{2}\wv^\top\xv_i\xv_i^\top\wv-\wv^\top\BB{b}
\end{align*}
 \deleted{Here, the condition numbers can be computed as $\kappa=\frac{\eig_{\mathrm{max}}\left(\Bm\right)}{\eig_\mathrm{min}\left(\Bm\right)}$ for standard first order methods and as $\kappa=\frac{\max_{i}{\Vnorm{\xv_i}^2}}{\eig_\mathrm{min}\left(\Bm\right)}$ for stochastic methods.} 
 \added{Thus, we can find an approximate solution to $\Bm^{-1}\BB{b}$ using either standard first-order methods or more advanced stochastic solvers. As we will see, using the advanced solvers together with good initializations will lead to a significant improvement in iteration complexity of our method.}

\subsection{\replaced{The top-k generalized eigenvector problem}{Problem statement}}
\begin{definition}[Top-$k$ generalized eigenvectors] Given symmetric matrices $\B{A}$ and $\B{B}$ where $\B{B}$ is positive definite, the top-$k$ generalized eigenvectors $\wv_1, \ldots,\wv_k$ are defined by
\begin{align}\label{def:ge}
\begin{aligned}
     \wv_i \in\, &\underset{\wv}{\argmax}
    & & 
    	\vert\wv^\top\B{A}\wv\vert
\\
    & \mbox{subject to}
    && \wv^\top\B{B}\wv=1
\\
    &&& \wv^\top\B{B}\wv_j=0 \quad \forall j \in \{1,2,\ldots,i-1\}.
\end{aligned}
\end{align}
\end{definition}
Problem~\eqref{def:ge} can be reduced to the problem of computing the top-$k$ eigenvectors associated with $k$ eigenvalues of largest magnitude of $\B{M}=\B{B}^{-1/2}\B{A}\B{B}^{-1/2}$ and then multiplying by $\B{B}^{-1/2}$ {(see, Lemma~2)}. Let us use the simple power method for this task. Let $\yv=\B{B}^{-1/2}\B{M}^T\B{x}$ be the result of $T$ iterations of the power method followed by multiplying {with} $\B{B}^{-1/2}$, we can write $\yv=\left(\B{B}^{-1}\B{A}\right)^T\B{B}^{-1/2}\xv$. Since the power method converges {from} a random initial vector $\xv$ with high probability, one can ignore $\B{B}^{-1/2}$ and compute instead $\yv=\left(\B{B}^{-1}\B{A}\right)^T\xv$. Thus, we can compute our desired eigenvectors by simply alternating between applying $\B{A}$ and $\B{B^{-1}}$ to a random initial vector. The authors of \cite{GJK16} suggested a procedure (GenELin) for this task which applies $\B{B}^{-1}$ approximately by least squares solvers. It is easy to see that GenLin simulates the noisy power method in \cite{HP14}, and hence it will converge as long as the noise is under control.
\section{Computing the Leading Generalized Eigenvector}
In this section, we introduce \texttt{NAPI}, an noisy accelerated power method for solving~\eqref{eq:EP}. We then characterize the convergence rate of the proposed algorithm for the special case of computing the leading generalized eigenvector.

\subsection{Description of the algorithm}
Our algorithm is designed to run in an inner-outer fashion and it hinges upon: (i) an acceleration technique \replaced{inspired by}{using} Chebyshev polynomial\added{s for approximating matrix powers}; and (ii) \replaced{an efficient inexact matrix-vector product computation}{a practical inexactness computation of matrix-vector products}. In particular, the acceleration technique in the outer loop\deleted{s} allows to reduce the number of outer iterations, while the subproblems in the inner loop\deleted{s can be solve} \added{are} efficiently \replaced{solved by some of the most advanced convex optimization solvers.}{by calling the most advanced solvers from convex optimization.} These solvers \replaced{are capable of exploiting sparsity and problem structure and scale to truly large-scale data sets.}{inherit  many advantages allowing to exploit sparsity in the data as well as the problem structures, and hence making our algorithm scalable.} \replaced{The last critical component of our algorithm is a well-designed warm start procedure which reduces the practical running times even further.}{Warm-start is also critical in our algorithm for further improve the running time of the algorithm.}

\begin{algorithm}[!t]
	\caption{Noisy Accelerated Power Iteration (\texttt{NAPI})}
	\begin{algorithmic}[1]\label{alg:1}
		\vspace{0.2cm}
		\REQUIRE Initial points $\wv_{-1}=0$ and $\wv_0$, and  parameter $\beta$.\\[0.0cm]
		\STATE $\wv_0 \gets \wv_0/\Mnorm{\wv_0}{\B{B}}$ \\[0.0cm]
		\FOR{$t=0,1,\ldots, T-1$}
		\STATE $\alpha_t \gets \wv_t^\top\B{A}\wv_t/\wv_t^\top\B{B}\wv_t$
		\STATE $\wtv_{t+1} \approx \argmin_{\wv \in \R^d}\left\{\frac{1}{2}\wv^\top \B{B}\wv - \wv^\top\B{A}\wv_t \right\} $ \\[0.1ex] \hspace{0.0cm}(initialize the solver with $\alpha_t  \wv_t$)\\[0.0cm]
		\STATE $\wtv_{t+1} \gets \wtv_{t+1} - \beta \wv_{t-1}$\\[0.0cm]
		\STATE $\wv_t \gets \frac{\wv_t}{\Mnorm{\wtv_{t+1}}{\B{B}}}$ and $\wv_{t+1} \gets \frac{\wtv_{t+1}}{\Mnorm{\wtv_{t+1}}{\B{B}}}$
		\ENDFOR
		\ENSURE $\wv_T$
	\end{algorithmic}
\end{algorithm}

\subsection{Deriving acceleration through Chebyshev magic}

As discussed above, \replaced{the power method can be accelerated by finding a good low-degree polynomial approximation to $x^T$.}{to improve on power method, it suffices to find a better low-degree approximation to the scalar function $x^T$.} However, \replaced{to evaluate a  polynomial of degree $T$, we may need to store up to $T$ terms before adding them up.}{it may require to store up to  $T$ terms at the $T$th iteration to make an approximation.} Fortunately, \replaced{the next result demonstrates that we can construct a Chebyshev polynomial recursively using only the two previous terms.}{this task can be done very efficiently by using only two previous terms thanks to the following extremely useful representation of the Chebyshev polynomial.}

\begin{definition}[Chebyshev Polynomials]\label{def:Chebyshev}
For a nonnegative integer $t$, the degree $t$ Chebyshev polynomial of the first kind\added{, $\Cheb{t}{z}$,} is defined recursively as follows:
\begin{align*}
	\Cheb{0}{z}&=1, 
	\,\,\,
	\Cheb{1}{z}=z,
	\\
	\Cheb{t}{z} &= 2z \Cheb{t-1}{z} - \Cheb{t-2}{z} \quad \text{for} \quad t \geq 2.
\end{align*} 
\end{definition}

To motivate our algorithm, let us focus on the unconstrained version of Problem~\eqref{def:ge}, where we relax the \added{condition that the solution should have} unit $\Bm$-norm\deleted{condition}.
Our goal is to find the leading eigenvector of $\Bm^{-1}\Am$ \replaced{by}{via} approximating $\left(\Bm^{-1}\Am\right)^T$ \replaced{using a}{by the} Chebyshev polynomial. Using Definition~\ref{def:Chebyshev}\deleted{of Chebyshev polynomial} with $z=\Bm^{-1}\Am$ and letting $\yv_t=\Cheb{t}{\Bm^{-1}\Am}\wv_0$, we can approximate $\left(\Bm^{-1}\Am\right)^T\wv_0$ by recursively \replaced{building}{build} a new iterate using only the previous two iterates:
\begin{align}\label{eq:unconstrained:recursive:y}
\yv_{t+1}=2\Bm^{-1}\Am\yv_t - \yv_{t-1}.
\end{align}
Note that the result in Proposition~\ref{prop:poly:existence} is stated for approximating over the interval $[-1,1]$. Therefore, we will consider a scaled Chebyshev polynomial of the form $C_t(\Bm^{-1}\Am)=\Cheb{t}{\frac{\Bm^{-1}\Am}{c}}$ for some constant $c$. Letting $\xv_t=C_t(\Bm^{-1}\Am)\wv_0$, we will thus study the recursion 
\begin{align}\label{eq:unconstrained:recursive}
\xv_{t+1}=\Bm^{-1}\Am\xv_t - \beta\xv_{t-1},
\end{align}
where $\beta$ is a constant depending on $c$. As we shall see, determining the best possible $\beta$ leads to an improved convergence rate for our method.

Returning to the constrained case, we wish to use a (scaled) Chebyshev polynomial to get an improved convergence rate while ensuring that the iterates have unit $\Bm$-norm. To this end, in our algorithm, we normalize the intermediate iterate $\wtv_{t+1}$ to have unit $\Bm$-norm and rescale the current iterate $\wv_t$ by the same factor as that for $\wtv_{t+1}$.  The iterates generated by  Algorithm~\ref{alg:1} can be written as
\begin{align}\label{eq:noisless:iterates}
&\wtv_{t+1} := \Bm^{-1}\Am\wv_t - \beta \wv_{t-1},\,\, \gamma_{t+1}=\Mnorm{\wtv_{t+1}}{\B{B}} \\ \label{eq:noiseless:scale}
&\wv_t := \wv_t \gamma_{t+1}^{-1}, \,\,
\wv_{t+1} := \wtv_{t+1} \gamma_{t+1}^{-1},
\end{align}
and the latest iterate $\wv_T$ can be expressed as 
\begin{align*}
\wv_T=\frac{C_T\left(\Bm^{-1}\Am\right)\wv_0}{\Mnorm{C_T\left(\Bm^{-1}\Am\right)\wv_0}{\Bm}},
\end{align*}
where $C_T$ is a scaled Chebyshev polynomial of degree $T$. That is, the output vector obtained by applying~\eqref{eq:noisless:iterates}-\eqref{eq:noiseless:scale} recursively is equivalent to performing only  recursion~\eqref{eq:unconstrained:recursive} without normalizing the iterates at each step, but only at the last step. To see that, suppose that $\wv_{-1}=\xv_{-1}=0$ and $\wv_0=\xv_0$, then we will show by induction that $\wv_t=\xv_t \alpha_t^{-1}$ for all $t\geq 1$, where $\alpha_t=\gamma_t \gamma_{t-1}\ldots\gamma_0 >0$. Assume that $\wv_i=\xv_i \alpha_i^{-1}$ for all $i \leq t$,   the sequence generated by~\eqref{eq:noisless:iterates}-\eqref{eq:noiseless:scale} can be written as  
\begin{align}
	\wv_{t+1} =\left( \Bm^{-1}\Am\wv_t - \beta \wv_{t-1} \gamma_t^{-1}\right)\gamma_{t+1}^{-1}.
\end{align}
It follows that 
\begin{align*}
	\wv_{t+1} &=\left( \Bm^{-1}\Am \xv_t \alpha_t^{-1} - \beta \xv_{t-1} \alpha_{t-1}^{-1} \gamma_t^{-1}\right)\gamma_{t+1}^{-1}
	\\
	&=
		\left( \Bm^{-1}\Am \xv_t \alpha_t^{-1} - \beta \xv_{t-1} \alpha_t^{-1} \right)\gamma_{t+1}^{-1}
	\\
	&=
		\left( \Bm^{-1}\Am \xv_t  - \beta \xv_{t-1}  \right)\alpha_t^{-1}\gamma_{t+1}^{-1}
	\\
	&=
	 \xv_{t+1}\alpha_{t+1}^{-1},
\end{align*}
as desired. The following result, whose proof can be found in~\cite[Theorem~8]{SHM17}, characterizes the convergence rate of Algorithm~\ref{alg:1} in the absence of noise.
\begin{proposition}[Convergence of the exact method]\label{prop:1}
Let $\eig_i$ be the eigenvalues of $\Bm^{-1}\Am$ satisfying $\norm{\eig_1} > \norm{\eig_2} \geq \ldots \geq \norm{\eig_d}$, and let $\gap=1-\norm{\eig_2}/\norm{\eig_1}$ be its relative eigenvalue gap. If the least squares subproblems in Step~4 of Algorithm~\ref{alg:1} are solved exactly and the parameter $\beta$ satisfies $\norm{\eig_2} \leq 2\sqrt{\beta} < \norm{\eig_1}$, then  
\begin{align*}
	\sin^2 \angl{}{\wv_T,\uv_1}
	\leq 
	4
	\left(
		\frac{
			2\sqrt{\beta}
		}{
			\norm{\eig_1} + \sqrt{\eig_1^2-4\beta}
		}
	\right)^{2T}\tan^2 \angl{}{\wv_{0},\uv_1}.
\end{align*}
Moreover, if $\beta=\eig_2^2/4$, then one can achieve an $\epsilon$-optimal solution, i.e, $\sin \angl{}{\wv_T,\uv_1} \leq \epsilon$, after at most $T = \mathcal{O}\left(\frac{1}{\sqrt{\gap}}\log\frac{\tan \theta_0}{\epsilon}\right)$ iterations.
\end{proposition}

The above result matches the iteration complexity of the Lanczos method, a significantly more complex algorithm \cite{GV12}. However, note that the result requires exact computation of $\Bm^{-1}\Am\wv$ for some vector $\wv$, which (as already discussed) is prohibitively expensive for large scale data sets. 

\subsection{Warm-start}
Initialization plays a critical role in designing an efficient method, cf. the work on the standard power iteration in~\cite{GJK16}. Recall that the matrix-vector product $\Bm^{-1}\Am\wv_t$ is computed approximately by solving the least squares problem
\begin{align*}
	\min_{\wv\in\R^d} f\left(\wv\right) := \frac{1}{2}\wv^\top \B{B}\wv - \wv^\top\B{A}\wv_t,
\end{align*} 
whose unique solution is $\wv^\star =\Bm^{-1}\Am\wv_t$. As the iterates converge, the solutions to the least squares problems in consecutive iterations become increasingly similar. It is natural to use  a scaled version of the previous iterate as initialization. This gives the following initial objective
\begin{align*}
	f\left(\alpha\wv_t\right)=\frac{1}{2}\alpha^2 \wv_t^\top\Bm\wv_t - \alpha\wv_t^\top\Am\wv_t.
\end{align*}
Minimizing $f\left(\alpha\wv_t\right)$ over $\alpha$ gives us the optimal scaling $\alpha_t^\star = \wv_t^\top\B{A}\wv_t/\wv_t^\top\B{B}\wv_t$. 
As will be shown later, with this initialization, we only need to run the least squares solver until the initial error has been reduced by a constant factor, independent of the final suboptimality required.
\subsection{Convergence argument}
In this subsection, we will show that as long as the errors in the matrix-vector   multiplications are controlled in an appropriate way, \texttt{NAPI} enjoys the same convergence rate as \replaced{the idealized method where subproblems are solved exactly. 
}{that of the noiseless case.}

Due to the \replaced{inexact solution of subproblems}{noises}, the output of Algorithm~\ref{alg:1} \replaced{is no longer equivalent to}{is no longer a result of} applying the Chebyshev polynomial to an initial vector, hence the proof of Proposition~\ref{prop:1} is not applicable {anymore}. However,\deleted{since} the iterates defined by Algorithm~\ref{alg:1} can be seen as a second order iterative system \replaced{and be}{with constant coefficient, it can be} cast into the \replaced{following form}{linear form as follows:}
\begin{align*}
	\augv_{t+1}:=
	\begin{bmatrix}
	\wv_{t+1}\\
	\wv_t
	\end{bmatrix}
	&= 		
	\frac{1}{\gamma_{t+1}}
	\left(
	\begin{bmatrix}
	\Bm^{-1}\Am & -\beta\IM{} \\[0.3em]
	\IM{} 		& \B{0} 
	\end{bmatrix}
	\begin{bmatrix}
	\wv_{t}\\
	\wv_{t-1}
	\end{bmatrix}
	+
	\begin{bmatrix}
	\nv_t\\
	\B{0}
	\end{bmatrix}
	\right)
	\\
	&= 
	\frac{1}{\gamma_{t+1}}
	\left(\Cm \augv_t + \anv_t\right),
\end{align*}
where \added[id=VM]{ $\nv_t=\wtv_{t+1}-\Bm^{-1}\Am\wv_t$ is the approximation error in Step~4 of Algorithm~1,}
\begin{align*}
\Cm
=
\begin{bmatrix}
\Bm^{-1}\Am & -\beta\IM{} \\[0.3em]
\IM{} 		& \B{0} 
\end{bmatrix}
\quad \text{and} \quad 
\anv_t
=
\begin{bmatrix}
\nv_t\\
\B{0}
\end{bmatrix}.
\end{align*}
\replaced{To establish the convergence rate of \texttt{NAPI}, we will thus study the behaviour of this extended system.}{Then, we can instead study the convergence behaviour of the extended system.} . {Note that the extended vectors $\augv_t$ do not have unit norm w.r.t any fixed matrix, hence care must be taken when performing the convergence analysis for this sequence.} We begin with the existence of eigenvectors and eigenvalues of $\Bm^{-1}\Am$ and $\Cm$.

\begin{lemma}[Existence of Eigenbasis]\label{lem:eig:base}
Let $\left(\BB{p}_i,\eig_i\right)$ be the eigenpairs of the symmetric matrix $\Bm^{-1/2}\Am\Bm^{-1/2}$. Then, $\left(\uv_i =\Bm^{-1/2}\BB{p}_i, \eig_i\right)$ is an eigenpair of $\Bm^{-1}\Am$. In addition,  
$\left(
	[\mu_i \uv_i^\top, \uv_i^\top]^\top
	, \mu_i
\right)
$  is an eigenpair of $\Cm$, where the eigenvalues $\mu_i$ satisfy    
\begin{align*}
	\mu_i^2 - \eig_i\mu_i + \beta = 0.
\end{align*}
\begin{proof}
See Appendix~\ref{appendix:A1}.
\end{proof}
\end{lemma}

The next lemma shows that under a certain condition\added{ on {the} suboptimality in solving the least-squares subproblems}, the angle between $\augv_t$ and the optimal eigenvector $\vv_1$ \replaced{shrinks geometrically in}{geometrically shrinks after} each iteration.
\begin{lemma}\label{lem:1}
Let $\vv_1$ be a leading eigenvector of the augmented matrix $\B{C}$ and let $\mu_i$ be the eigenvalues of $\Cm$ such that $\norm{\mu_1} > \norm{\mu_2} \geq \ldots \ge \norm{\mu_{2d}}$. If the errors in solving the least squares problems satisfy
\begin{align}\label{eq:noise:condition}
\hspace{-0.2cm}	\Vnorm{\nv_t}_{\Bm}  \leq \frac{\norm{\mu_1}-\norm{\mu_2}}{4}\min\left\{\sin\angl{}{\augv_t,\vv_1}, \cos\angl{}{\augv_t,\vv_1}\right\}
\end{align}
for all $t \in \N$, then 
\begin{align*}
	\tan\angl{}{\augv_{t+1},\vv_1} 
	\leq
		\frac{\norm{\mu_1}+3\norm{\mu_2}}{3\norm{\mu_1}+ \norm{\mu_2}}
		\tan\angl{}{\augv_t,\vv_1}. 
\end{align*}
\begin{proof}
See Appendix~\ref{appendix:A2}.
\end{proof}
\end{lemma}

We are now ready to state our main theorem\replaced{, quantifying}{that shows} the iteration complexity of \replaced{\texttt{NAPI}:}{
the proposed algorithm.}
\begin{theorem}\label{thrm:1}
 Let $\eig_i$ be the eigenvalues of $\Bm^{-1}\Am$ satisfying $\norm{\eig_1} > \norm{\eig_2} \geq \ldots \geq \norm{\eig_d}$, and let $\gap=1-\norm{\eig_2}/\norm{\eig_1}$ be its relative eigenvalue gap. For $\norm{\eig_2} \leq 2\sqrt{\beta} < \norm{\eig_1}$, the output of Algorithm~\ref{alg:1} satisfies
\begin{align*}
	\sin^2 \angl{}{\augv_T,\vv_1}
	\leq 
	\left(
		\frac{1}{2}
		+
		\frac{
			\sqrt{\beta}
		}{
			\norm{\eig_1} + \sqrt{\eig_1^2-4\beta}
		}
	\right)^{2T}\tan^2 \theta_0,
\end{align*}
where $\theta_0=\angl{}{\augv_0,\vv_1}$. Furthermore, if the parameter $\beta$ is set to $\norm{\eig_2}^2/4$, we can achieve an $\epsilon$-optimal solution after at most
\begin{align*}
	T = \mathcal{O}\left(\frac{1}{\sqrt{\gap}}\log\frac{\tan \theta_0}{\epsilon}\right)
\end{align*}  iterations.  
The total running time of the algorithm  is 
\begin{align*}
\mathcal{O}\left(
\frac{1}{\sqrt{\gap}}
\left(
\mathcal{T}\left(\frac{\gap \cos^2 \theta_0}{32}\right)
\log\left(\tan \theta_0\right)
+
\mathcal{T}\left(\frac{\gap}{32}\right)
\log \frac{1}{\epsilon}
\right)
+
\frac{\mathrm{nnz}\left(\Am\right)+\mathrm{nnz}\left(\Bm\right)}{\sqrt{\gap}}
\log\frac{\tan \theta_0}{\epsilon}
\right),
\end{align*}
where $\mathcal{T}\left(\tau\right)$ is the time that a least squares solver takes to reduce the residual error by a factor $\tau$.
\begin{proof}
See Appendix~\ref{appendix:B}
\end{proof}
 \end{theorem}

\begin{remark}\label{rmk:1}
The result in Theorem~\ref{thrm:1} for the extended system readily implies the convergence of $\wv_t$ to an optimal solution to Problem~\ref{def:ge}. More specifically,\deleted{by} invoking Lemma~\ref{lem:angle:space} with 
\begin{align*}
	{\B{X} = \BB{v}_1/\Vnorm{\BB{v}_1}_{\bigB},}	
	\B{Y}=\augv_T/\Mnorm{\augv_T}{\bigB},
	\,\, \mbox{and} \,\,
	\Bm
	=
		\bigB
	= 
		\rm{\begin{bmatrix}
			\Bm & \B{0} \\[0.0em]
			\B{0} 		& \Bm 
		\end{bmatrix}},\,\,	
\end{align*}
yields
\begin{align*}
\min_{\alpha\in\{\pm 1\}}
\Vnorm{
	\frac{\bigB^{\frac{1}{2}}}{\sqrt{\gamma_T^{-2}+1}}
	\begin{bmatrix}
	\wv_T\\
	\wv_{T-1}
	\end{bmatrix}
	-
	\frac{\alpha\bigB^{\frac{1}{2}}}{\sqrt{\mu_1^2+1}}
	\begin{bmatrix}
	\mu_1 \uv_1\\
	\uv_1
	\end{bmatrix}
}
\leq
2\sin\frac{\angl{}{\augv_T,\vv_1}}{2}.
\end{align*}
It follows that
\begin{align*}
\min_{\alpha\in\{\pm 1\}}
	\Vnorm{
   		\frac{\gamma_T}{\sqrt{\gamma_T^2+1}}
		\Bm^{\frac{1}{2}}
		\wv_T
		-
		\frac{\alpha\mu_1}{\sqrt{\mu_1^2+1}}
		\Bm^{\frac{1}{2}}
		\uv_1
	}
	\leq
		2\sin\frac{\angl{}{\augv_T,\vv_1}}{2}.
\end{align*}
{Since $\Vnorm{\wv_T}_\Bm=\Vnorm{\uv_1}_\Bm=1$, the above inequality indicates that $\wv_T$ converges to a leading eigenvector ($\uv_1$ or $-\uv_1$) of $\Bm^{-1}\Am$, an optimal solution of Problem~\ref{def:ge}.}

\end{remark}

\begin{remark}
Theorem~\ref{thrm:1} shows that as long as the \replaced{errors}{noises} arising from  solving the least squares problems \added{approximately} are controlled \replaced{appropriately}{in an appropriate way}, the proposed algorithm enjoys the same convergence rate as\deleted{that of} the \replaced{idealized method}{noiseless case}(c.f. Proposition~\ref{prop:1}). Oddly, our analysis shows that the accelerated method is less sensitive to errors than the standard power method in~\cite{GJK16}. \replaced{In particular}{Particularly}, the condition~\eqref{eq:noise:condition} in Lemma~\ref{lem:1} requires the noises \added{to} scale with $\sqrt{\gap}$ while the condition in~\cite{GJK16} requires the noises scale with $\gap$. {It is worth mentioning that compared to~\cite{GJK16}, our method has one additional parameter (the momentum gain), whose optimal value is not known in general. This situation is quite common in momentum-like methods such as the heavy-ball method~\cite{Pol87} in convex optimization. However, it has been widely observed that even with a suboptimal estimate of the momentum parameter, one often can achieve significantly speed-ups over the original algorithm.} 
\end{remark}

\begin{corollary}
If we use Nesterov's accelerated method as a least squares solver, the running time of Algorithm~1 to achieve an $\epsilon$-optimal solution is bounded by
\begin{align*}
	\mathcal{O}\left(
		\frac{\mathrm{nnz}\left(\Bm\right) \sqrt{\kappa}}{\sqrt{\gap}}
		m_1
		+
		m_2
	\right),
\end{align*}
where 
\begin{align*}
m_1&=\left(
\log\frac{1}{\cos \theta_0}\log\frac{1}{\gap\cos \theta_0}
+
\log\frac{1}{\epsilon}
\log\frac{1}{\gap}
\right)
\\
m_2&=\frac{\mathrm{nnz}\left(\Am\right)}{\sqrt{\gap}}\log\frac{1}{\epsilon\cos \theta_0}.
\end{align*}
This is a linear running time in the input size of the problem, i.e., the number of nonzero elements of $\Am$ and $\Bm$.
\end{corollary}


\subsection{Extension to the Top-$k$ Generalized Eigenvectors}
We consider an extension of our accelerated method \replaced{to the problem of}{for} finding the top-$k$ generalized eigenvectors. The algorithm is formally given as Algorithm~\ref{alg:2}. It is a natural generalization of Algorithm~\ref{alg:1}, where the iterates are now matrices of size $d \times k$, and the normalization steps are replaced by orthogonalization steps using the Gram-Schmidt procedure with the inner product $\InP{\cdot}{\cdot}_{\Bm}$. Similarly as in the top-$1$ case,  we will study the convergence behaviour of the sequence generated by Algorithm~\ref{alg:2} via\deleted{studying} the extended system given by
\begin{align*}
	\Ym_{t+1}:=
	\begin{bmatrix}
	\Wm_{t+1}\\
	\Wm_t
	\end{bmatrix}
	= 		
	\left(
	\begin{bmatrix}
	\Bm^{-1}\Am & -\beta\IM{} \\[0.3em]
	\IM{} 		& \B{0} 
	\end{bmatrix}
	\begin{bmatrix}
	\Wm_{t}\\
	\Wm_{t-1}
	\end{bmatrix}
	+
	\begin{bmatrix}
	\nvk_{t}\\
	\B{0}
	\end{bmatrix}
	\right)
	\Rm_{t+1}^{-1}.
\end{align*}

\begin{algorithm}[!t]
	\caption{\texttt{NAPI} for top-$k$ Generalized Eigenvectors}
	\begin{algorithmic}[1]\label{alg:2}
		\vspace{0.2cm}
		\REQUIRE Initial points $\Wm_{-1}=\B{0}$, $k$, $T$, and  $\beta$.\\[0.0cm]
		\STATE Let $\Wm_0$ be a random orthonormal basis w.r.t $\InP{\cdot}{\cdot}_{\Bm}$\\[0.0cm]
		\FOR{$t=0,1,\ldots, T-1$}
		\STATE $\Zm_t \gets \left(\Wm_t^\top\B{B}\Wm_t\right)^{-1}\Wm_t^\top\B{A}\Wm_t$
		\STATE $\Wtm_{t+1} \approx \argmin_{\Wm }\left\{\frac{1}{2}\tr\left(\Wm^\top \B{B}\Wm\right) - \tr\left(\Wm^\top\B{A}\Wm_t \right)\right\}$ (initialize the solver with $\Wm_t\Zm_t $)\\[0.0cm]
		\STATE $\Wtm_{t+1} \gets \Wtm_{t+1} - \beta \Wm_{t-1}$\\[0.0cm]
		\STATE $ \Wm_{t+1} \gets \Wtm_{t+1}\Rm_{t+1}^{-1}$ via QR-factorization w.r.t $\InP{\cdot}{\cdot}_{\Bm}$
		\STATE $\Wm_t \gets \Wm_t \Rm_{t+1}^{-1}$
		\ENDFOR
		\ENSURE $\Wm_{s}$
	\end{algorithmic}
\end{algorithm}

{Let $\Um=[\uv_1,\ldots,\uv_k]$ with $\Um^\top\Bm\Um=\IM{}$ be the matrix of top-$k$  eigenvectors of $\Bm^{-1}\Am$ and let $\eig_1,\ldots,\eig_k$ be the corresponding top-$k$ eigenvalues. Then, by Lemma~\ref{lem:eig:base}, the matrices of top-$k$ eigenvectors $\bar{\Vm}$ and eigenvalues $\Eeig$ of the extended matrix $\Cm$ are $\bar{\Vm}=[\Eeig\Um^\top, \Um^\top]^\top$ and $\Eeig = \diag{\mu_1,\ldots,\mu_k}$ with $\mu_i$ being solutions to the equations $\mu_i^2 - \eig_i\mu_i + \beta = 0$, $i=1\ldots,k$. Let $\Vm=\bar{\Vm}(\IM{}+ \Eeig^2)^{-1/2}$,  and let $$\bigB
=
\begin{bmatrix}
\Bm & \B{0} \\[0.0em]
\B{0} 		& \Bm 
\end{bmatrix},
$$
then $\Vm^\top\bigB\,\Vm=\IM{k}$. Furthermore, using QR decomposition w.r.t $\bigB$, one can write $\Ym_t$ as $\Ym_t=\bar{\Ym}_t\B{G}_t$ for some matrix $\bar{\Ym}_t$ satisfying $\bar{\Ym}_t^\top\bigB\,\bar{\Ym}_t=\IM{}$ and for some nonsingular matrix $\B{G}_t\in \R^{k\times k}$.}

The following theorem characterizes the convergence properties of the above extended system. Note that the angles in the theorem are measured w.r.t $\bigB$. 

\begin{theorem}\label{thrm:2}
Let $\eig_i$ be the eigenvalues of $\Bm^{-1}\Am$ satisfying $\norm{\eig_1} \geq \ldots \geq \norm{\eig_k} > \norm{\eig_{t+1}} \geq \ldots \geq \norm{\eig_d}$, and let $\gap_k=1-\norm{\eig_{k+1}}/\norm{\eig_k}$ be its relative eigenvalue gap. 
If the errors in solving the least squares problems satisfy
\begin{align*}
	\Vnorm{\B{E}_t}_\Bm
	\leq \frac{\norm{\mu_k}-\norm{\mu_{k+1}}}{4}\min\left\{\sin\theta({\bar{\Ym}}_t,\Vm),\cos\theta({\bar{\Ym}}_t,\Vm)\right\},
\end{align*}
for all $t \in \N$, and if $\norm{\eig_{k+1}} \leq 2\sqrt{\beta} < \norm{\eig_k}$, the output of Algorithm~\ref{alg:2} satisfies
\begin{align*}
	\sin^2 \angl{}{\bar{\Ym}_{T},\Vm}
	\leq
	\left(
	\frac{1}{2}
	+
	\frac{
		\sqrt{\beta}
	}{
		\norm{\eig_k} + \sqrt{\eig_k^2-4\beta}
	}
	\right)^{2T}\tan^2 \theta_0,
\end{align*}
where $\theta_0=\theta(\bar{\Ym}_0,\Vm)$.
Furthermore, if the parameter $\beta$ is set to $\norm{\eig_{k+1}}^2/4$, we can achieve an $\epsilon$-optimal solution after at most
\begin{align*}
	T = \mathcal{O}\left(\frac{1}{\sqrt{\gap_k}}\log\frac{\tan \theta_0}{\epsilon}\right)
	\end{align*}  iterations.  
	The total running time of the algorithm  is 
	\begin{align*}
	&\mathcal{O}\left(
	\frac{1}{\sqrt{\gap_k}}
	\left(
	\mathcal{T}\left(\frac{\gap_k \cos^4 \theta_0}{128k \gamma^2}\right)
	\log\left(\tan \theta_0\right)
	+
	\mathcal{T}\left(\frac{\gap_k}{128\gamma^2}\right)
	\log \frac{1}{\epsilon}
	\right)
	\right.
	\\
	&\hspace{0.5cm}
	\left.
	+
	\frac{1}{\sqrt{\gap_k}}\left(
	\mathrm{nnz}\left(\Am\right)k
	+
	\mathrm{nnz}\left(\Bm\right)k
	+
	dk^2
	\right)
	\log\frac{\tan \theta_0}{\epsilon}
	\right),
	\end{align*}
	where $\gamma=\eig_1/\eig_k$ and  $\mathcal{T}\left(\tau\right)$ is the time a least squares solver takes to reduce the residual error by a factor $\tau$.
\begin{proof}
See Appendix~\ref{appendix:C}.
\end{proof}
\end{theorem}

{For $k=1$, we have shown in Remark~\ref{rmk:1} that the convergence of the extended sequence translates to the convergence of the original one.  However, for a general case, we need a more involved analysis to derive a similar conclusion. We start by invoking Lemma~\ref{lem:angle:space} with 
$\B{X}=\bar{\Ym}_T$, $\B{Y}=\Vm$, $\Bm=\bigB$ to obtain
\begin{align*}
\min_{ \B{Q}^\top\B{Q}=\IM{k}}
\Vnorm{
	\bigB^{\frac{1}{2}}
	\begin{bmatrix}
	\Wm_T\\
	\Wm_{T-1}
	\end{bmatrix}
	\B{G}_T
	-	
	\bigB^{\frac{1}{2}}
	\begin{bmatrix}
	\Um \Eeig\\
	\Um
	\end{bmatrix}
	(\IM{}+ \Eeig^2)^{-1/2}\B{Q}
}
	\leq
		2\sin\frac{\angl{}{\bar{\Ym}_{T},\Vm}}{2}.
\end{align*}
It follows that there exists an orthogonal matrix $\B{Q}$ such that
\begin{align*}
\Vnorm{
	\Wm_T
	\B{G}_T
	-		
	\Um \Eeig
	(\IM{}+ \Eeig^2)^{-1/2}\B{Q}
	}_{\Bm}
	\leq
		2\sin\frac{\angl{}{\bar{\Ym}_{T},\Vm}}{2}.
\end{align*}
Therefore, if we let $\B{M}=\Eeig(\IM{}+ \Eeig^2)^{-1/2}\B{Q}$, then $\Wm_T$ can be written as
\begin{align*}
	\Wm_T
	=
		\Um \Eeig
		(\IM{}+ \Eeig^2)^{-1/2}\B{Q}
		\B{G}_T^{-1}
		+
		\B{E}\B{G}_T^{-1}
	=
		\big(\Um+\B{E}\B{M}^{-1}\big)\B{M}\B{G}_T^{-1},		
\end{align*}
where $\B{E}$ is a perturbation matrix satisfying $\Vnorm{\B{E}}_\Bm\leq 2\sin\frac{\angl{}{\bar{\Ym}_{T},\Vm}}{2}$.
Note that since $\B{M}\B{G}_T^{-1}$ is nonsingular, the range space of $\Wm_T$ is identical to that of $\Um+\B{E}\B{M}^{-1}$. Consider the following matrix 
\begin{align*}
	\B{P}= 
		\big(\Um+\B{E}\B{M}^{-1}\big)
		\big[
			\big(\Um+\B{E}\B{M}^{-1}\big)^\top
			\Bm
			\big(\Um+\B{E}\B{M}^{-1}\big)
		\big]^{-1/2},
\end{align*}
then it can be verified that $\B{P}^\top\Bm\B{P}=\IM{}$, thereby being an orthogonal basis for the range space of $\Wm_T$. Using the binomial expansion identity~\cite{Hig08}
\begin{align*}
	(\IM{}-\B{A})^{-1/2}
	=
		\IM{}
		+\frac{1}{2}\B{A}
		+
		\rm{H.O.T}
		,
\end{align*}
for $\Vnorm{\B{A}}<1$ and the fact that $\Um^\top\Bm\Um=\IM{}$, one can write the matrix $\B{P}$ as $\B{P}=\big(\Um+\B{E}\B{M}^{-1}\big)\big(\IM{}+\B{E}'\big)$ with $\Vnorm{\B{E}'}_{\Bm}=O(\Vnorm{\B{E}}_\Bm)$. Therefore, if we denote by $\Um_\perp$ an orthogonal basis (w.r.t $\Bm$) of the subspace spanned by $\uv_{k+1},\ldots,\uv_{d}$, then it holds that
\begin{align*}
	\Um_\perp^\top\Bm\B{P} 
	= 
		\Um_\perp^\top\Bm\B{E}\B{M}^{-1}
		+
		\Um_\perp^\top\Bm\B{E}\B{M}^{-1}\B{E}',
\end{align*}
which implies that
\begin{align*}
	\sin\theta(\Wm_T,\Um) 
	= 
		\Vnorm{\Um_\perp^\top\Bm\B{P}}
	= O(\Vnorm{\B{E}}_{\Bm})
	\leq c \sin\frac{\angl{}{\bar{\Ym}_{T},\Vm}}{2},
\end{align*}
for some constant $c$. The last inequality indicates that the original sequence converges at the same rate as the extended one, as desired.}

\section{Application to CCA}
Our main result can readily lead to improvement of several downstream applications. In this section, we will discuss the benefits of our proposed method in solving {the} CCA problem.

In CCA, we are provided with pairs of data points from two views $(\xv_1,\yv_1), \ldots
, (\xv_n,\yv_n)$ where $\xv_i \in \R^{1\times d_1}$,  $\yv_i \in \R^{1\times d_2}$, and $n$ is the size of the training set. For example, the views can be the visual image data and audio data in a video of people speaking. Let  $\B{X} = \left[\xv_1, \ldots , \xv_n\right]^\top \in \R^{n \times d_1}$ and $\B{Y} = \left[\yv_1, \ldots , \yv_n\right]^\top \in \R^{n \times d_2}$ be the data matrices for each view, respectively. We define the relevant empirical covariance matrices as $\Cov{11}=\frac{1}{n}\B{X}^\top \B{X} + \reg_1\IM{}$, $\Cov{22}=\frac{1}{n}\B{Y}^\top\B{Y}+ \reg_2\IM{}$, and  $\Cov{12}=\frac{1}{n}\B{X}^\top\B{Y}$, where $\reg_1$ and $\reg_2$ are regularization parameters. 
The goal is to find a linear transformation of the points in each view that retains as much as the redundant information between the views. More specifically, CCA finds the solution of the following problem
\begin{align}\label{P:CCA}
\begin{aligned}
    \left(\cX_i, \cY_i\right) &\in \underset{\cX, \cY}{\argmax}
    & & 
    \hspace{-0.2cm}	\cX^\top\Cov{12}\cY
\\
    &\hspace{-0.2cm}\text{subject to}
    &&\hspace{-0.2cm} \cX^\top\Cov{11}
    \cX= \cY^\top\Cov{22}\cY=1
\\
	&&& \hspace{-0.2cm}\cX^\top\Cov{11}\cX_j=\cY^\top\Cov{22}\cY_j=0, \,\, \forall j \leq i-1. 
\end{aligned}
\end{align}
It is easy to check that any stationary point $\left(\cX_i, \cY_i\right)$ of Problem~\eqref{P:CCA} satisfies 
\begin{align}\label{CCA:sys}
	\Cov{12}\cY_i = \eig_i \Cov{11} \cX_i \quad \text{and} \quad 
	\Cov{12}^\top\cX_i = \eig'_i \Cov{22} \cY_i,
\end{align}
where $\eig_i$ and $\eig'_i$ are two constants. Using the constraints, it can be \replaced{verified}{checked} that $\eig_i=\eig'_i$. The system~\eqref{CCA:sys} above can \added{therefore} be cast into the following generalized eigenproblem 
\begin{align}\label{eq:redu:cca}
	\begin{bmatrix}
		\B{0} & \Cov{12} \\[0.0em]
		\Cov{12}^\top 		& \B{0} 
	\end{bmatrix}
	\begin{bmatrix}
    	\cX_i\\
   		\cY_i
	\end{bmatrix}
	=
	\eig_i
	\begin{bmatrix}
		 \Cov{11} & \B{0}\\[0.0em]
		 \B{0} 	& \Cov{22}
	\end{bmatrix}
	\begin{bmatrix}
    	\cX_i\\
   		\cY_i
	\end{bmatrix}.
\end{align}
The solution to this problem has $d_1+d_2$ eigenvalues $\eig_1 > \eig_2 > \ldots > -\eig_2 > -\eig_1$; for each eigenvalue $\eig_i$ with the corresponding eigenvector 
$[ 	\cX_i^\top,	\cY_i^\top]^\top$, $-\eig_i$ is also an eigenvalue with the corresponding eigenvector 
$[\cX_i^\top, 	-\cY_i^\top]^\top$.
Let $\Am$ and $\Bm$ be the matrices on the left and right in~\eqref{eq:redu:cca}, respectively, then the eigenvalue gap of $\Bm^{-1}\Am$ is zero. Thus, iterative methods for finding the leading eigenvector of $\Bm^{-1}\Am$ may not converge. Since the top $2k$-dimensional eigen-space of $\Bm^{-1}\Am$ is the subspace spanned by $[ 	\cX_i^\top,	\cY_i^\top]^\top$ and 
$[\cX_i^\top, 	-\cY_i^\top]^\top$, $i=1,\ldots,k$, one can solve the above problem for the top-$2k$ eigenvector problem and then {choose} any orthogonal basis that spans the output subspace and a random $k$-dimensional projection of those vectors. Finally, to fulfil the constraints in CCA problem, we perform the last orthogonalization steps w.r.t $\Cov{11}$ and $\Cov{22}$, respectively. The formal description of the discussions above are given in Algorithm~\ref{alg:3}.
\begin{algorithm}[!h]
\caption{Noisy Accelerated Power Method for CCA}
 \begin{algorithmic}[1]\label{alg:3}
 \vspace{0.2cm}
\REQUIRE Data matrices $\B{X}, \B{Y}$, $T$, $k$, and $\beta$.\\[0.1cm]
\STATE $ \B{A} \gets 
	\begin{bmatrix}
			\B{0} & \Cov{12} \\[0.0em]
			\Cov{12}^\top 		& \B{0} 
\end{bmatrix}, 
\quad
	\B{B} \gets 
	\begin{bmatrix}
			\Cov{11} & \B{0} \\[0.0em]
			 \B{0}   & \Cov{22}
\end{bmatrix}
$\\[0.1cm]

\STATE $[\tilde{\CX}_s^\top,\tilde{\CY}_s^\top]^\top \gets \texttt{NAPI}\left(\Am,\Bm, \beta, T, 2k\right)$
\STATE Generate a $2k\times k$ random Gaussian matrix $\BB{G}$	
\STATE $\CX_T \gets \tilde{\CX}_T \BB{G}$ and  $\CY_T \gets \tilde{\CY}_T \BB{G}$ 
\STATE $\CX_T \gets \CX_T \B{R}_{\mathrm{x}}^{-1}$ via QR-factorazation w.r.t $\InP{\cdot}{\cdot}_{\Cov{11}}$
\STATE $\CY_T \gets \CY_T \B{R}_{\mathrm{y}}^{-1}$ via QR-factorazation w.r.t $\InP{\cdot}{\cdot}_{\Cov{22}}$
\ENSURE $\CX_T$, $\CY_T$
\end{algorithmic}
\end{algorithm}

We remark that when running Algorithm~\ref{alg:3}, we do not need to form $\Am$ and $\Bm$ explicitly thanks to their diagonal forms. In fact, we only need to apply $\Cov{11}$, $\Cov{12}$, and $\Cov{22}$ as operators, which can be done efficiently by applying $\B{X}$ and $\B{Y}$ appropriately. For example, at some points of the algorithm, one needs to compute $\xv^\top\Cov{11}\yv$ for some vectors $\xv,\yv$, one can compute $\B{X}\yv$ and $\B{X}\xv$ and then compute the inner product of those vectors without forming $\Cov{11}$. This also allows to exploit the sparsity of the data matrices $\B{X}$ and $\B{Y}$, which is critical for large-scale {problems}. 

The least squares problem at the $t$th iteration of procedure $\texttt{NAPI}$ can be written explicitly as follows
\begin{align*}
\min_{\CX,\CY}  f_t\left(\CX,\CY\right)
&:=
\frac{1}{2}
\tr\left(
\begin{bmatrix}
\CX\\
\CY
\end{bmatrix}^\top
\begin{bmatrix}
\Cov{11} & \B{0} \\[0.0em]
\B{0}   & \Cov{22}
\end{bmatrix}
\begin{bmatrix}
\CX\\
\CY
\end{bmatrix}
\right)
-
\tr\left(
\CX^\top
\Cov{12}	
\CX_{k-1}
-
\CY^\top
\Cov{12}^\top	
\CY_{k-1}\right).
\end{align*}
The \emph{condition number} of this problem is $$\kappa = \max\left(\frac{\eig_{\mathrm{max}}\left(\Cov{11}\right)}{\eig_{\mathrm{min}}\left(\Cov{11}\right)}, \frac{\eig_{\mathrm{max}}\left(\Cov{22}\right)}{\eig_{\mathrm{min}}\left(\Cov{22}\right)}\right).$$
If one wishes to use stochastic optimization solvers, one can be express the objective function above as 
\begin{align*}
	\min_{\CX,\CY} f_t\left(\CX,\CY\right) := \frac{1}{n} \sum_{i=1}^{n}f_t^i\left(\CX,\CY\right),
\end{align*}
where
\begin{align*}
f_t^i\left(\CX,\CY\right)
&=
\frac{1}{2}
\tr\left(
\begin{bmatrix}
\CX\\
\CY
\end{bmatrix}^\top
\begin{bmatrix}
\xv_i^\top\xv_i + \reg_1 \IM{} & \B{0} \\[0.0em]
\B{0}   & \yv_i^\top\yv_i + \reg_2 \IM{}
\end{bmatrix}
\begin{bmatrix}
\CX\\
\CY
\end{bmatrix}
\right)
-
\tr\left(
\CX^\top
\Cov{12}	
\CX_{k-1}
-
\CY^\top
\Cov{12}^\top	
\CY_{k-1}\right).
\end{align*}
In this case, the relevant condition number depends on the individual \emph{Lipschitz} constant of the component functions, and can be computed as 
$$\tilde{\kappa} = \max\left(\frac{\max_i\left(\Vnorm{\xv_i}^2\right)}{\eig_{\mathrm{min}}\left(\Cov{11}\right)}, \frac{\max_i\left(\Vnorm{\yv_i}^2\right)}{\eig_{\mathrm{min}}\left(\Cov{22}\right)}\right).$$

The following lemma shows that Steps~5 and 6 in Algorithm~\ref{alg:3} do not cause much loss in the alignment with the true canonical space. 
\begin{lemma}{\cite[Lemma~14]{GJK16}}\label{lem:Ge}
If the output of procedure~\texttt{NAPI} has the angle at most $\theta$	with the true top-$2k$ generalized eigenspace of $\Bm$, $\Am$, then with probability at least $1-\delta$, both $\CX_T$ and $\CY_T$ has angle at most $\mathcal{O}\left(k^2\theta/\delta^2\right)$ with the true top-$k$ canonical space of $\B{X}$ and $\B{Y}$.
\end{lemma}

Finally, we state the following theorem charactering the iteration complexity of Algorithm~\ref{alg:3}.
\begin{theorem}\label{thrm:3}
Given data matrices $\B{X}$ and $\B{Y}$, let $\Am$ and $\Bm$ be the matrices defined in Algorithm~\ref{alg:3}. Let $\eig_i$ be the eigenvalues of $\Bm^{-1}\Am$ satisfying $\eig_1 > \eig_2 > \ldots > -\eig_2 > -\eig_1$, and let $\gap_k=1-{\norm{\eig_{k+1}}}/{\norm{\eig_k}}$ be its relative eigenvalue gap. Let $\CX^\star$ and $\CY^\star$ be the top-$k$ true canonical space of $\B{X}$ and $\B{Y}$, respectively. If $\beta = \norm{\eig_{k+1}}^2/4$, then with probability at least $1-\delta$, the output of Algorithm~\ref{alg:3} satisfying
\begin{align*}
	\min\left\{
		\sin \angl{}{\CX_T,\CX^\star}, 
		\sin \angl{}{\CY_T,\CY^\star}
	\right\}
	\leq \epsilon,
\end{align*} 
in time
\begin{align*}
&\mathcal{O}\left(
\frac{1}{\sqrt{\gap_k}}
\mathcal{T}\left(\frac{c\gap_k \delta^2 \cos^4 \theta_0}{k^3 \gamma^2}\right)
\log\left(\tan \theta_0\right)			
+
\frac{1}{\sqrt{\gap_k}}
\mathcal{T}\left(\frac{c\gap_k \delta^2}{k^2\gamma^2}\right)
\log \frac{1}{\epsilon}
+
\frac{
	\mathrm{nnz}\left(\B{X},\B{Y}\right) k + dk^2
}{
	\sqrt{\gap_k}
}
\log\frac{\tan \theta_0}{\epsilon}
\right),
\end{align*}
where $\mathrm{nnz}\left(\B{X},\B{Y}\right)=\mathrm{nnz}\left(\B{X}\right)+\mathrm{nnz}\left(\B{Y}\right)$, $\gamma=\eig_1/\eig_k$ and $c$ is a universal constant.
\begin{proof}
	The proof follows readily from Lemma~\ref{lem:Ge} and Theorem~\ref{thrm:2}.
\end{proof}
\end{theorem}

\begin{corollary}
	If we use Nesterov's accelerated method as a least squares solver, the running time of Algorithm~3 to achieve an $\epsilon$-optimal solution is bounded by
	\begin{align*}
	\mathcal{O}\left(
	\frac{\mathrm{nnz}\left(\B{X},\B{Y}\right) k\sqrt{\kappa}}{\sqrt{\gap_k}}
	m_1
	+
	m_2
	\right).
	\end{align*}
	On the other hand, if we choose to use SVRG, then the running time is bounded by
	\begin{align*}
	\mathcal{O}\left(
	\frac{k\,\mathrm{nnz}\left(\B{X},\B{Y}\right) \left( 1 + \sqrt{\tilde{\kappa}/{n}}\right)}{\sqrt{\gap_k}}
	m_1
	+
	m_2
	\right),
	\end{align*}
	where 
\begin{align*}
m_1&=\left(
\log\frac{1}{\cos \theta_0}\log\frac{k\gamma}{\delta \gap_k \cos \theta_0}
+
\log\frac{1}{\epsilon}
\log\frac{k\gamma}{\delta\gap_k}
\right)
\\
m_2&=\frac{dk^2}{\sqrt{\gap_k}}\log\frac{1}{\epsilon\cos \theta_0}.
\end{align*}	
	
\end{corollary}

\section{Experimental Results}
In this section, we perform experiments on real-word data sets to validate the effectiveness of the proposed algorithms. In the plots, we illustrate the performance of \texttt{NAPI}-based CCA on the following real-world data sets.

{\bf{Mediamill}} is an annotated video data set from the Mediamill Challenge for automated detection of semantic concepts, which contains 85 hours of international broadcast news data \cite{SWV06}. These news videos are first automatically segmented into 30,993 subshots. Each image is a representative key-frame of a video subshot annotated with 101 labels and consists of 120 features. CCA is used to learn the correlation structure between the images and its labels.

{\bf{MNIST}} is a data set of 60,000 handwritten digits from 0 to 9 \cite{LBY98}. Each digit is represented by an image of $392 \times 392$ values in $\left[0,1\right]$. CCA is performed to explore the correlated representations between the left and right halves of the handwritten digit images.

{\bf{XRMB}} is Wisconsin Xray Microbeam database consits of simultaneous acoustic and articulatory recordings \cite{WMW90, WAL15}.  The inputs to CCA are the acoustic and articulatory features concatenated over a 7-frame window around each frame, giving acoustic vectors $\BB{x}\in \R^{273}$ and articulatory vectors $\BB{y}\in \R^{112}$.
\begin{table*}[t]
    \caption{
                Brief Summary of Data sets
    }
    \label{table:1}
\begin{center}
    \begin{tabular}{c|c|c|c|c}
        \hline \hline
        Data set & Description & $d_1$ & $d_2$ & $n$
         \bigstrut\\
        \hline
        Mediamill &   Image and its labels  &   101 & 120 & 30,993 \bigstrut\\
        XRMB &   Acoustic and articulation measurements   &   273 & 112 & 1,429,236  \\[0.3ex]
        MNIST &  Left and right halves of images    &   392 & 392 &  60,000\\[0.3ex]
        \hline
    \end{tabular}
\end{center}
\end{table*}

\begin{figure*}[t!]
	\centering \subfigure[Mediamill, $\reg_1=\reg_2=10^{-3}$]{
		{\includegraphics[width=0.32\textwidth]{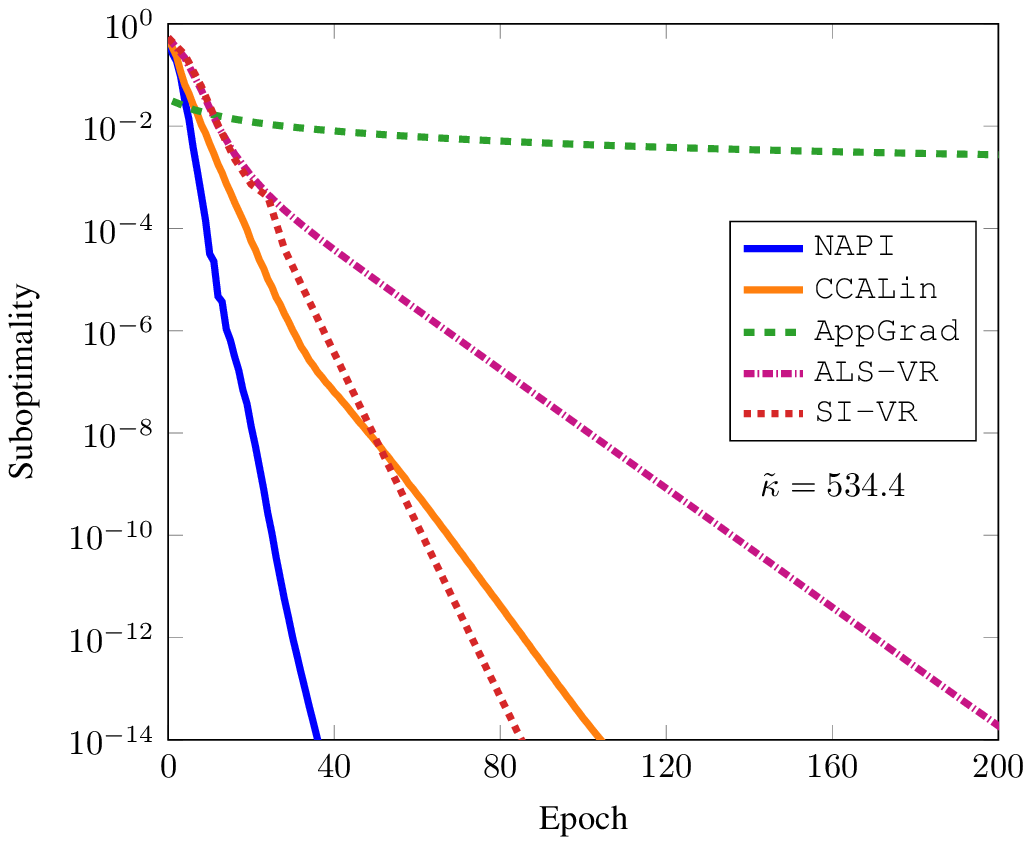}}
		\label{fig:1a} }~\subfigure[MNIST, $\reg_1=\reg_2=10^{-3}$]{
		{\includegraphics[width=0.32\textwidth]{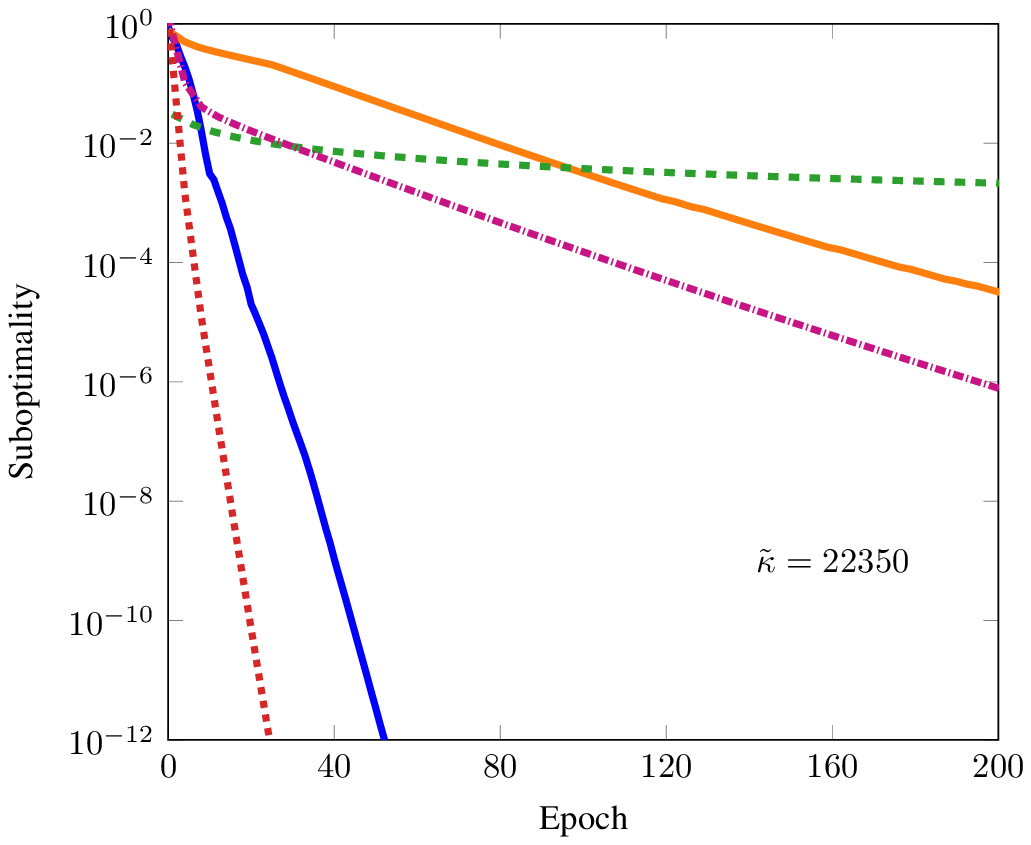}}
		\label{fig:1b} }~ \subfigure[XRMB, $\reg_1=\reg_2=10^{-3}$]{
		{\includegraphics[width=0.32\textwidth]{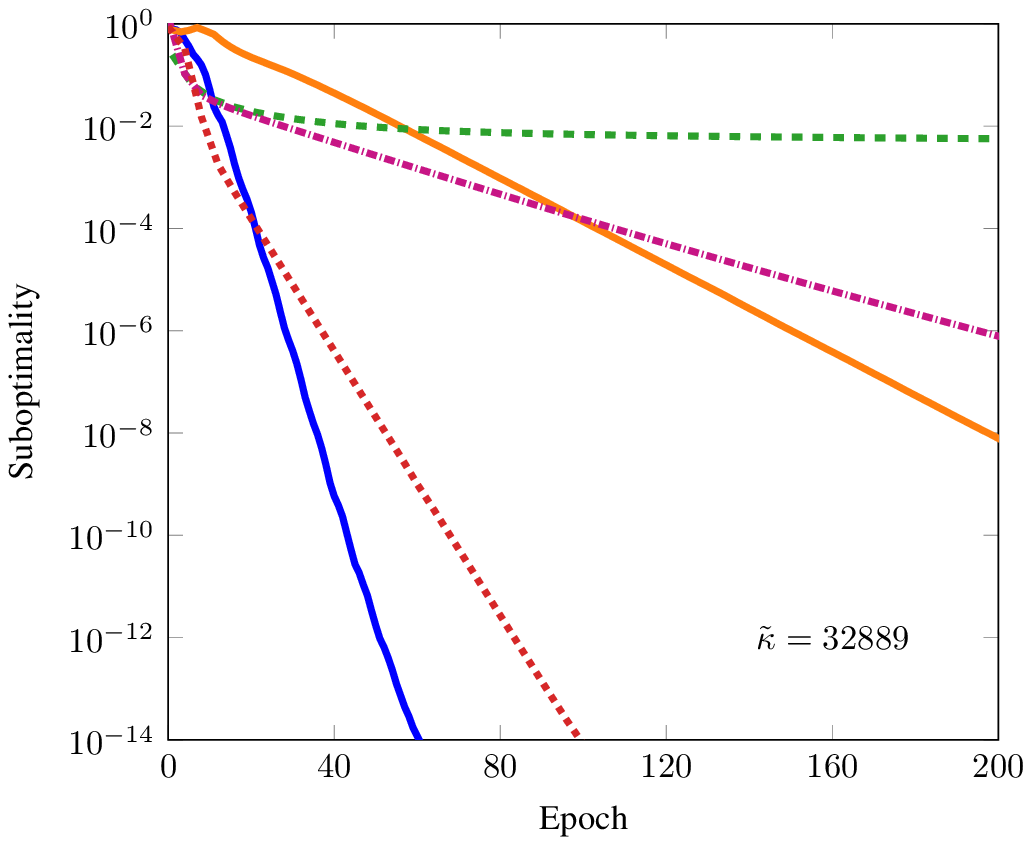}}
		\label{fig:c} }\par\subfigure[Mediamill, $\reg_1=\reg_2=10^{-5}$]{
		{\includegraphics[width=0.32\textwidth]{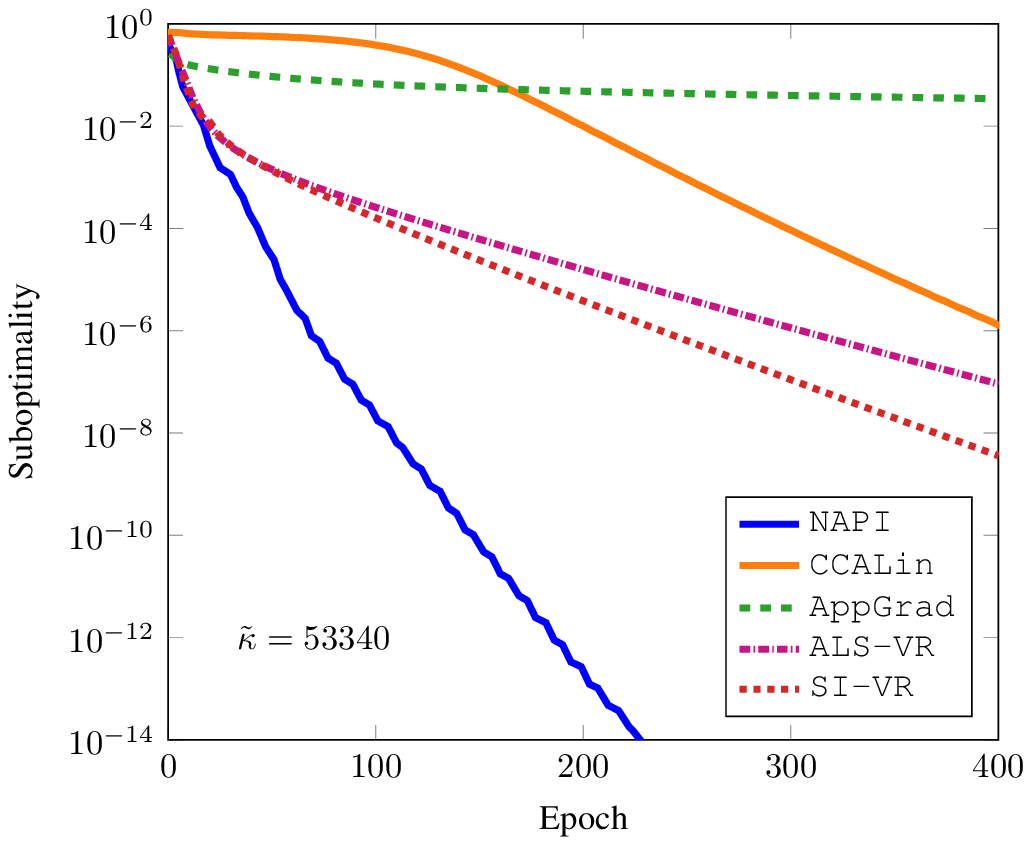}}
		\label{fig:2a} }~\subfigure[MNIST, $\reg_1=\reg_2=10^{-5}$]{
		{\includegraphics[width=0.32\textwidth]{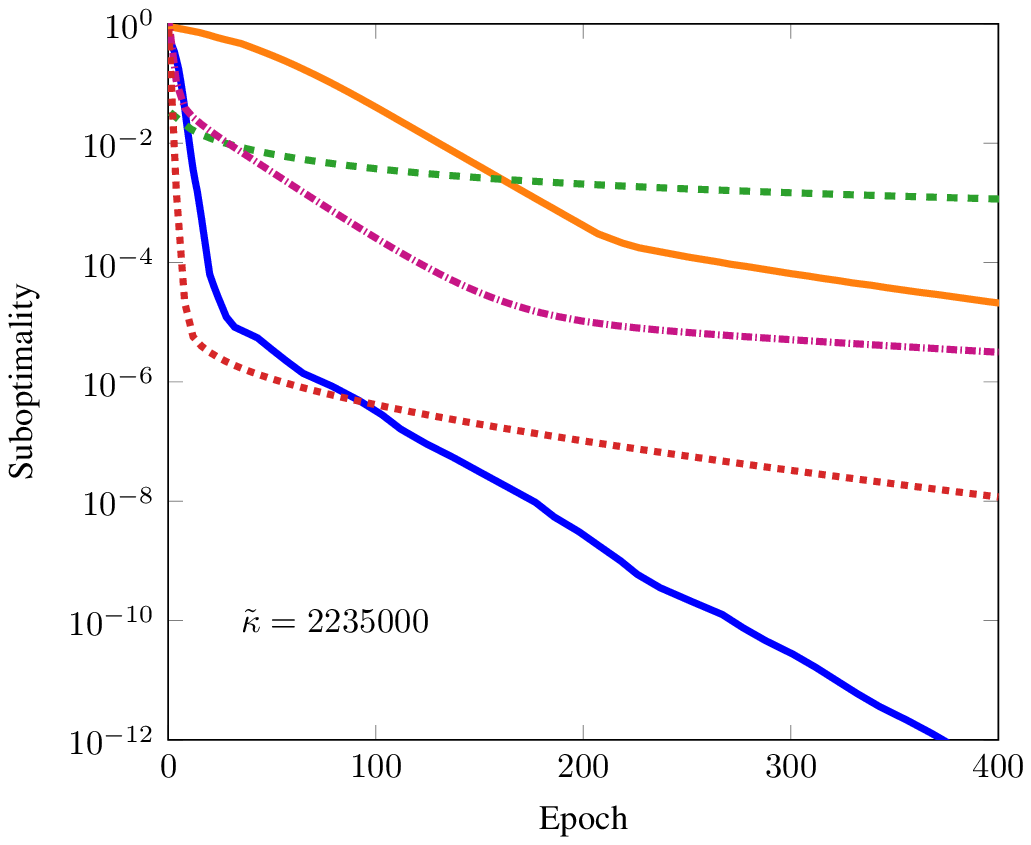}}
		\label{fig:2b} }~\subfigure[XRMB, $\reg_1=\reg_2=10^{-5}$]{
		{\includegraphics[width=0.32\textwidth]{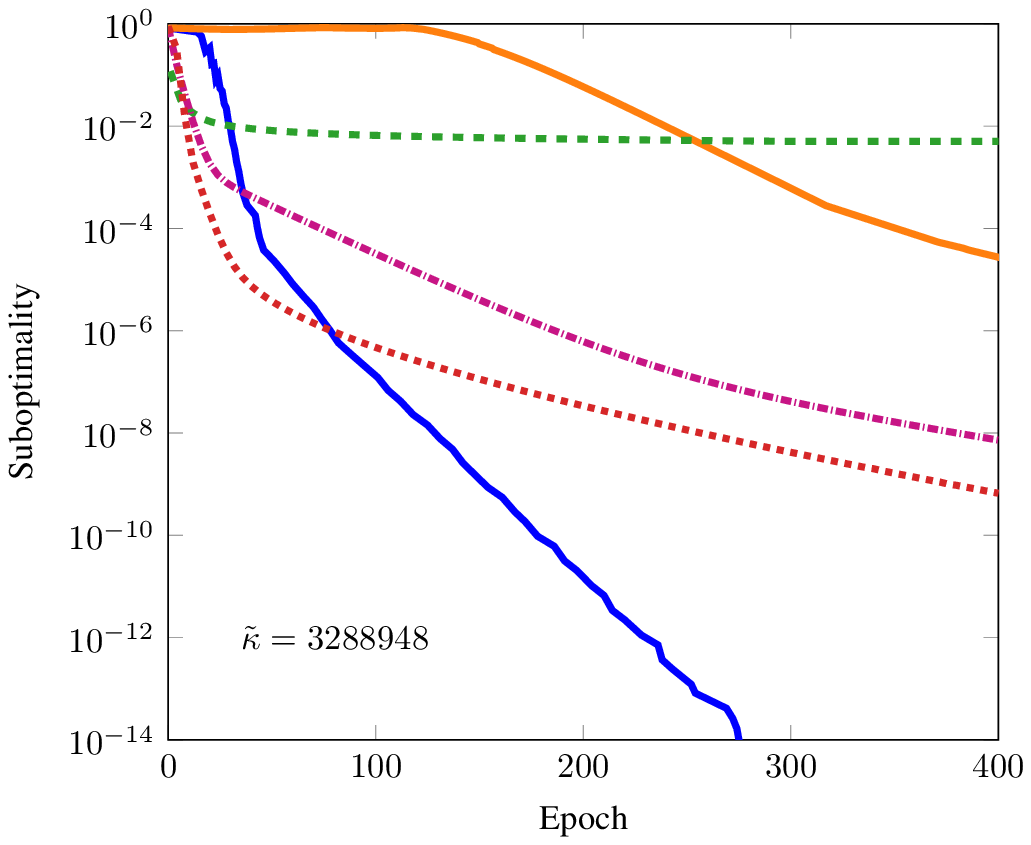}}
		\label{fig:2c} }~\caption{Suboptimality versus the number of iterations for different algorithms solving the CCA problem when $k=1$, and the regularization parameters $\reg_1$ and $\reg_2$ are given in the plot.} \label{fig:3}
\end{figure*}

Figure~\ref{fig:3} shows the suboptimality in objective versus the number of epochs (passes over the full data set) for different algorithms solving CCA problem with $k=1$ starting from a random initial vector. We compare the following algorithms: \texttt{AppGrad} \cite{MLF15}, \texttt{ALS-VR} \cite{WWG16}, \texttt{CCALin} \cite{GJK16}, {and \texttt{SI-VR} \cite{WWG16}} where SVRG is used as the least squares solver. For \texttt{CCALin}, and \texttt{NAPI}, we set the SVRG  stepsizes to be $1/\max_{ i}\left(\Vnorm{\xv_i}^2,\Vnorm{\yv_i}^2\right)$. {We use the implementation of  \texttt{SI-VR} provided by the authors in \cite{WWG16}. This implementation uses slightly different parameters than analyzed in \cite{WWG16}, but outperforms the theoretically suggested values.} For \texttt{AppGrad}, we tune the stepsizes and report the result from using the best ones.
For each dataset, we vary the regularization parameters to influence the condition numbers of the least squares problems, where larger regularization implies better conditioning. 

{It can be seen that our proposed method converges linearly to an optimal solution and significantly outperforms the other algorithms, especially for ill-conditioned problems. As predicted by the theoretical results, \texttt{NAPI} and \texttt{SI-VR} have superior performance
compared to other non-accelerated methods. It should be emphasized that for well-conditioned problems such as the case of the \texttt{Mediamill} data set  with $\reg_x=\reg_y=10^{-3}$, \texttt{CCALin} and \texttt{AlS-VR} can have comparable performance to \texttt{NAPI} and \texttt{SI-VR}. Finally, we observe that for ill-conditioned problems, the convergence of \texttt{SI-VR} can be quite slow with a similar slope as \texttt{ALS-VR}, while \texttt{NAPI} still maintains an accelerated rate and quickly finds a high accuracy solution.}


\section{Conclusion}
We have proposed and analyzed a simple algorithm for finding the dominant generalized eigenvectors of a pair of symmetric matrices. The algorithm admits a linear convergence rate and was shown to outperform the current state-of-the-art algorithms. To achieve that result, our algorithm exploits the problem structure and decomposes the overall problem into a sequence of subproblems. The approximation theory is used to accelerate the overall process which helps to achieve the asymptotic iteration complexity of the Lanczos method, while each of the subproblem can be solved efficiently by advanced iterative solvers and with a well-design initialization. This translates into a \emph{linear} time algorithm in the input sizes, suitable for many large-scale applications in machine learning and data science.

\section*{Acknowledgment}
This research was sponsored in part by the Knut and Alice Wallenberg Foundation and the Swedish Research Council. We thank Weiran Wang for sharing with us the implementation of the shift-and-invert method in~\cite{WWG16} and two anonymous reviewers for their useful comments and suggestions.

\bibliographystyle{plain}
\bibliography{refs}

\appendix

\section{Proof of Auxiliary Lemmas}

\subsection{Proof of Lemma~\ref{lem:eig:base}}\label{appendix:A1}
For the first part, we have 
\begin{align*}
\Bm^{-1}\Am \left(\Bm^{-1/2}\BB{p}_i\right)
	=\Bm^{-1/2}\left(\Bm^{-1/2}\Am \Bm^{-1/2}\BB{p}_i\right)
	=\eig_i \Bm^{-1/2}\BB{p}_i.
\end{align*}
For the second part, for any eigenpair $\left(\uv,\eig\right)$ of $\Bm^{-1}\Am$, let $\mu$ be a solution of $\mu^2 - \eig\mu + \beta = 0$. Then, 
\begin{align*}
	\Cm\begin{bmatrix}
	\mu \uv\\
	\uv
	\end{bmatrix}
	= 
	\begin{bmatrix}
	\Bm^{-1}\Am & -\beta\IM{} \\[0.3em]
	\IM{} 		& \B{0} 
	\end{bmatrix}
	\begin{bmatrix}
	\mu \uv\\
	\uv
	\end{bmatrix}	
	=
	\begin{bmatrix}
	\added{(}\mu\Bm^{-1}\Am  -\beta\IM{}\added{)}\uv \\[0.3em]
	\mu \uv
	\end{bmatrix}
	= 
	\mu \begin{bmatrix}
	\mu \uv\\
	\uv
	\end{bmatrix}.
\end{align*}
Therefore, the vector $\vv=[\mu \uv^\top, \uv^\top]^\top$ is an eigenvector of $\Cm$ \added{with eigenvalue $\mu$}. The proof is complete.
\subsection{Proof of Lemma~\ref{lem:1}}\label{appendix:A2}

{The proof follows closely those in~\cite{GJK16,HP14} with a few modifications to take the additional normalization steps in Algorithm~\ref{alg:1} into account.}
\replaced{First}{Firstly}, we recall that 
\begin{align*}
	\augv_{t+1}
	= 		
	\left(
	\Cm
	\augv_t
	+
	\anv_t
	\right)\gamma_{t+1}^{-1},
\end{align*}
where 
$
		\augv_t
		=
		[\wv_t^\top ,	\wv_{t-1}^\top]^\top
$, $\anv_t=[\nv_t^\top, \B{0}]^\top$, and $\gamma_{t+1}$ is a normalization term.		
{Let} $${\vv_1 = \frac{1}{\sqrt{\mu_1^2+1}}\left[\, \mu_1 \uv_1^\top, \uv_1^\top \right]^\top \,\,\, \mbox{and}\,\,\,\bigB
=
\begin{bmatrix}
\Bm & \B{0} \\[0.0em]
\B{0} 		& \Bm 
\end{bmatrix},}
$$
then $\Mnorm{\vv_1}{\bigB}=1$ since $\Mnorm{\uv_1}{\Bm}=1$. \added[id=VM]{Note that $\bigB$ is positive definite, hence the $\bigB$-norm is well defined.} Let $\Pv$ be a projection onto the space spanned by the top eigenvector of $\Cm$ and let $\PPv$ be a projection  onto the space spanned by the remaining ones, \added[id=VM]{with respect to the $\bigB$-norm}. It can be verified that $\Pv=\vv_1\vv_1^\top\bigB$ and $\PPv=\IM{}-\vv_1\vv_1^\top\bigB$. By \added{the} definition of principal angles under \added{the} $\mathcal{B}$-norm \added[id=VM]{for $k=1$, and the facts that $\Pv=\Pv^2$ and $\PPv=(\PPv)^2$, it is readily verified that}
$\cos  \angl{}{\augv_{t},\vv_1}=\frac{\Mnorm{\Pv\augv_t}{\bigB}}{\Mnorm{\augv_t}{\bigB}}$ 
and 
$\sin \angl{}{\augv_{t},\vv_1}=\frac{\Mnorm{\PPv\augv_t}{\bigB}}{\Mnorm{\augv_t}{\bigB}}$.

We now start by expanding $\tan\angl{}{\augv_{t+1},\vv_1}$ as follows
\begin{align}\label{eq:proof:lem:1:1}
\tan\angl{}{\augv_{t+1},\vv_1} 
&=
\frac{\Mnorm{\PPv \augv_{t+1}}{\bigB}}{\Mnorm{\Pv \augv_{t+1}}{\bigB}}
\nonumber\\
&=
\frac{\Mnorm{\PPv\left(\Cm\augv_t + \anv_t\right)}{\bigB}}{\Mnorm{\Pv \left(\Cm\augv_t + \anv_t\right)}{\bigB}}
\nonumber\\
&\leq
\frac{\Mnorm{\PPv\Cm\augv_t}{\bigB} + \Mnorm{\PPv\anv_t}{\bigB}}{\Mnorm{\Pv\Cm\augv_t}{\bigB} - \Mnorm{\Pv\anv_t}{\bigB}}.
\end{align}
Since
\begin{align*}
\Mnorm{\Pv\Cm\augv_t}{\bigB}
&\geq
\norm{\mu_1}\Mnorm{\Pv\augv_t}{\bigB}
\\
\Mnorm{\PPv\Cm\augv_t}{\bigB}
&=
\sqrt{\augv_t^\top\Cm^\top (\PPv)^\top\bigB \PPv \Cm \augv_t}
\leq
\norm{\mu_2}\Mnorm{\PPv\augv_t}{\bigB},
\end{align*} 
it follows that 
\begin{align}\label{eq:proof:lem:1:2}
\tan\angl{}{\augv_{t+1},\vv_1} 
&=
\frac{\Mnorm{\PPv\augv_t}{\bigB}}{\Mnorm{\Pv\augv_t}{\bigB}}
\times
\frac{
	\norm{\mu_2}
	+
	\frac{
		\Mnorm{\PPv\anv_t}{\bigB}
	}{
		\Mnorm{\PPv\augv_t}{\bigB}
	}
}{
	\norm{\mu_1}
	-
	\frac{
		\Mnorm{\Pv\anv_t}{\bigB} 
	}{
		\Mnorm{\Pv\augv_t}{\bigB}
	}
}
\nonumber\\
&\hspace{0.3cm}\leq
\tan\angl{}{\augv_{t},\vv_1} 
\times
\frac{
	\norm{\mu_2}
	+
	\frac{
		\Mnorm{\PPv\anv_t}{\bigB}
	}{
		\Mnorm{\PPv\augv_t}{\bigB}
	}
}{
	\norm{\mu_1}
	-
	\frac{
		\Mnorm{\Pv\anv_t}{\bigB} 
	}{
		\Mnorm{\Pv\augv_t}{\bigB}
	}
}
\nonumber\\
&\hspace{0.3cm}=
\tan\angl{}{\augv_{t},\vv_1} 
\times
\frac{
	\norm{\mu_2}
	+
	\frac{
		\Mnorm{\anv_t}{\bigB}
	}{
		\Mnorm{\augv_t}{\bigB}\sin\angl{}{\augv_{t},\vv_1}
	}
}{
	\norm{\mu_1}
	-
	\frac{
		\Mnorm{\anv_t}{\bigB} 
	}{
		\Mnorm{\augv_t}{\bigB}
		\cos \angl{}{\augv_{t},\vv_1}
	}
}.
\end{align}
To get a geometric rate of convergence, it suffices to choose the noise level such that the rightmost term in~\eqref{eq:proof:lem:1:2} is less than 1. Since $\Mnorm{\augv_t}{\bigB}=\sqrt{1+\gamma_t^{-2}} >1$, choosing 
\begin{align*}
\Vnorm{\nv_t}_{\Bm}  \leq \frac{\norm{\mu_1}-\norm{\mu_2}}{4}\min\left\{\sin\angl{}{\augv_t,\vv_1}, \cos\angl{}{\augv_t,\vv_1}\right\},
\end{align*}
yields
\begin{align}\label{eq:proof:lem:1:3}
\tan\angl{}{\augv_{t+1},\vv_1} 
\leq
\frac{\norm{\mu_1}+3\norm{\mu_2}}{3\norm{\mu_1}+\norm{\mu_2}}
\times
\tan\angl{}{\augv_{t},\vv_1},
\end{align}
completing the proof of Lemma~\ref{lem:1}.	

\section{Proof of Theorem~\ref{thrm:1}}\label{appendix:B}	
First, by applying  inequality~\eqref{eq:proof:lem:1:3} recursively,
we obtain
\begin{align}\label{eq:proof:thrm:1:rec}
\tan\angl{}{\augv_{t},\vv_1} 
&\leq
\left(\frac{\norm{\mu_1}+3\norm{\mu_2}}{3\norm{\mu_1}+\norm{\mu_2}}\right)^t
\times
\tan\angl{}{\augv_0,\vv_1}
\nonumber\\
&\leq
\left(1-\frac{1}{2}\left(1-\frac{\norm{\mu_2}}{\norm{\mu_1}}\right)\right)^t
\times 
\tan\angl{}{\augv_0,\vv_1},
\end{align}
\added[id=VM]{where the last step follows since $\norm{\mu_1}>\norm{\mu_2}$.}
We next examine the relative eigenvalue gap of the extended matrix $\Cm$.
For a fixed $i$, from Lemma~\ref{lem:eig:base}, the eigenvalues $\mu_i$ of $\Cm$ are the solutions of the characteristic equation
\begin{align*}
\mu_i^2 - \eig_i\mu_i + \beta = 0.
\end{align*}
If $2\sqrt{\beta} > \norm{\eig_1}$,  the roots of the characteristic equation are imaginary, and both have magnitude $2\sqrt{\beta}$, which results in a zero eigenvalue gap. On the other hand, if $2\sqrt{\beta} < \norm{\eig_2}$, then \deleted{the eigenvalue gap becomes}
\begin{align*}
1-\frac{\norm{\mu_2}}{\norm{\mu_1}}
&=
1 
-
\frac
{
	\norm{\eig_2} + \sqrt{\eig_2^2-4\beta}
}{
	\norm{\eig_1} + \sqrt{\eig_1^2-4\beta}
}
\end{align*} 
\replaced{while}{Finally}, for $\norm{\eig_2} \leq 2\sqrt{\beta} < \norm{\eig_1}$, we get

\begin{align}\label{eq:proof:thrm:1:gap}
1-\frac{\norm{\mu_2}}{\norm{\mu_1}}
=
1
-
\frac{
	2\sqrt{\beta}
}{
	\norm{\eig_1} + \sqrt{\eig_1^2-4\beta}
}.
\end{align}
We focus on the last \replaced{case, since this is the  range for $\beta$ which can offer accleration}{possibility since we may get an improved convergence rate from that choice of $\beta$}. 
\added[id=VM]{Note that the right-hand side of~\eqref{eq:proof:thrm:1:gap} is a decreasing function of $\beta$, hence achieving its optimal solution at $\beta=\norm{\eig_2}^2/4$.}
Combining~\eqref{eq:proof:thrm:1:rec} and~\eqref{eq:proof:thrm:1:gap} yields 
\begin{align}\label{eq:proof:thrm:1:2}
\sin^2 \angl{}{\augv_{t},\vv_1} 
\leq 
\tan^2 \angl{}{\augv_{t},\vv_1}
\leq 
\left(
\frac{1}{2}
+
\frac{
	\sqrt{\beta}
}{
	\norm{\eig_1} + \sqrt{\eig_1^2-4\beta}
}
\right)^{2t}
\tan^2\angl{}{\augv_0,\vv_1}.
\end{align}

We now pay attention to the running time of Algorithm~\ref{alg:1} required to achieve an $\epsilon$-optimal solution, when \replaced{the parameter $\beta$ is chosen optimally.}{the optimal parameter $\beta$ is chosen.} For $\beta=\norm{\eig_2}^2/4$, we have 
\begin{align}\label{eq:ext:gap}
1-\frac{\norm{\mu_2}}{\norm{\mu_1}}
&=
1
-
\frac{
	2\sqrt{\beta}
}{
	\norm{\eig_1} + \sqrt{\eig_1^2-4\beta}
}
\nonumber\\
&=
\frac{
	\norm{\eig_1}-\norm{\eig_2} + \sqrt{\eig_1^2-\eig_2^2}
}{
	\norm{\eig_1} + \sqrt{\eig_1^2-\eig_2^2}
}
\nonumber\\
&=
\frac{
	\sqrt{\norm{\eig_1}-\norm{\eig_2}}
	\left(
	\sqrt{\norm{\eig_1}+\norm{\eig_2}}
	-
	\sqrt{\norm{\eig_1}-\norm{\eig_2}}
	\right)
}{
	\norm{\eig_2}
}
\nonumber\\
&=		
\frac{
	2\sqrt{\gap}\sqrt{\norm{\eig_1}}
}{		
	\sqrt{\norm{\eig_1}+\norm{\eig_2}}
	+
	\sqrt{\norm{\eig_1}-\norm{\eig_2}}
}
\geq
\frac{\sqrt{\gap}}{\sqrt{2}}.
\end{align}
From~\eqref{eq:proof:thrm:1:rec}, it suffices to choose $k$ such that
\begin{align}\label{eq:iter:complx}
t \geq \frac{2\norm{\mu_1}}{\norm{\mu_1}-\norm{\mu_2}}\log\frac{\tan\angl{}{\augv_0,\vv_1}}{\epsilon}
\geq
\frac{2\sqrt{2}}{\sqrt{\gap}} \log\frac{\tan\angl{}{\augv_0,\vv_1}}{\epsilon}.
\end{align}

We now turn to control the errors arising from solving the least squares problems inexactly. Define the error in solving the least squares problem as 
\begin{align*}
\textsf{r}\left(\wv\right) = \Mnorm{\wv-\Bm^{-1}\Am\wv_t}{\Bm}^2=2\left(f\left(\wv\right)-f\left(\Bm^{-1}\Am\wv_t\right)\right).
\end{align*}
Let $\ires=\textsf{r}\left(\alpha_t\wv_t\right)$ and $\dres=\Mnorm{\nv_t}{\Bm}^2=\textsf{r}\left(\wtv_{t+1}\right)$ be the initial and required residual errors in Step~4 of Algorithm~\ref{alg:1}, respectively. Then, with our choice of initialization, we have by~\cite{GJK16} that $\ires\leq\eig_1^2\sin^2\angl{}{\wv_t,\uv_1}$.
Therefore, to guarantee the noise condition in~\eqref{eq:noise:condition}, it suffices to enforce that
\begin{align*}
\dres \leq \frac{\left(\norm{\mu_1}-\norm{\mu_2}\right)^2}{16}\min\left\{\sin^2\angl{}{\augv_t,\vv_1}, \cos^2\angl{}{\augv_t,\vv_1}\right\}.
\end{align*}
Thus, we only need to solve the least squares problem until the ratio of the final to initial error satisfies 
\begin{align}\label{eq:res:condition}
\frac{\dres}{\ires}
\leq
\frac{
	\left(\norm{\mu_1}-\norm{\mu_2}\right)^2
}{
	16 \eig_1^2
}
\frac{
	\min\left\{\sin^2\angl{}{\augv_t,\vv_1}, \cos^2\angl{}{\augv_t,\vv_1}\right\}
}{
	\sin^2\angl{}{\wv_t,\uv_1}
}.
\end{align}
For $\norm{\eig_2}\leq 2\sqrt{\beta} < \norm{\eig_1}$, we have $\norm{\mu_1} \geq \norm{\eig_1}$, hence if the errors above satisfy 
\begin{align}\label{eq:res:condition:new}
\frac{\dres}{\ires}
&\leq
\left(
1
-			
\frac{
	\norm{\mu_2}
}{
	\norm{\mu_1} 
}
\right)^2
\frac{\min\left\{\sin^2\angl{}{\augv_t,\vv_1}, \cos^2\angl{}{\augv_t,\vv_1}\right\}}{16 \sin^2\angl{}{\wv_t,\uv_1}}
\nonumber\\
&=
\frac{
	\gap
}{
	32 
}	
\min\left\{\frac{\sin^2\angl{}{\augv_t,\vv_1}}{\sin^2\angl{}{\wv_t,\uv_1}}, \frac{\cos^2\angl{}{\augv_t,\vv_1}}{\sin^2\angl{}{\wv_t,\uv_1}}\right\},
\end{align}
the condition~\eqref{eq:res:condition} is automatically satisfied. 

\added[id=VM]{We can decompose the running time into two phases: i) \emph{initial phase} due to large initial angle; and ii) \emph{convergence phase} due to high accuracy $\epsilon$ required.
In the initial phase, the angle are large, the min in~\eqref{eq:res:condition:new} is determined by the second term. Thus, we can set the ratio to be $\frac{\gap \cos^2\angl{}{\augv_0,\vv_1}}{32}$, until $\tan\angl{}{\augv_t,\vv_1}$ reduces to 1. This phase takes $\mathcal{O}\left(\frac{\log \tan\angl{}{\augv_0,\vv_1}}{\sqrt{\gap}}\right)$ iterations corresponding to $\epsilon=1$ in~\eqref{eq:iter:complx}. In the convergence phase, the angles are small, the first term determines the min in~\eqref{eq:res:condition:new}, we can set the ratio to be $\frac{\gap}{32}$ until we reach the target accuracy. This phase takes 
$\mathcal{O}\left(\frac{1}{\sqrt{\gap}} \log \frac{1}{\epsilon}\right)$ iterations corresponding to $\tan\angl{}{\augv_0,\vv_1}=1$ in~\eqref{eq:res:condition:new}.}

{Therefore, if we let $\mathcal{T}\left(\tau\right)$ be the time that a least squares solver takes to reduce the residual error by a factor $\tau$, then the total running time of Step~4 in Algorithm~\ref{alg:1} is given by 
\begin{align*}
\mathcal{O}\left(
\frac{1}{\sqrt{\gap}}
\left(
\mathcal{T}\left(\frac{\gap \cos^2 \theta_0}{32}\right)
\log\left(\tan \theta_0\right)
+
\mathcal{T}\left(\frac{\gap}{32}\right)
\log \frac{1}{\epsilon}
\right)
\right).
\end{align*}
Finally, the running time for Steps~3 and~6 is  
\begin{align*}
\mathcal{O}\left(
\frac{
	\mathrm{nnz}\left(\Am\right)
	+
	\mathrm{nnz}\left(\Bm\right)
}{
	\sqrt{\gap}
}
\log\frac{\tan \theta_0}{\epsilon}
\right),
\end{align*}
which completes the proof of Theorem~\ref{thrm:1}.}

\section{{Proof of Theorem~\ref{thrm:2}}}\label{appendix:C}	
To begin with, recall that $\Um$ is the matrix of top-$k$ eigenvectors of $\Bm^{-1}\Am$ and $\bar{\Vm}=[\Eeig\Um^\top, \Um^\top]^\top$ is the top-$k$ eigenvectors of the extended matrix $\Cm$ with  $\Eeig = \diag{\mu_1,\ldots,\mu_k}$ being the corresponding matrix of eigenvalues. Let $$\Vm=\bar{\Vm}(\IM{}+ \Eeig^2)^{-1/2}\,\,\, \mbox{and}\,\,\,
\bigB
=
\begin{bmatrix}
\Bm & \B{0} \\[0.0em]
\B{0} 		& \Bm 
\end{bmatrix},
$$
then $\Vm^\top\bigB\,\Vm=\IM{k}$. If we further let $\Vmp$ be an orthogonal basis w.r.t $\bigB$ of the orthogonal complement of $\rm{span}(\Vm)$, then one can decompose $\Cm$ as 
\begin{align}\label{eq:C:decomp}
\Cm=\Vm\Eeig\Vm^\top\bigB+\Vmp\Eeigp\Vmp^\top\bigB,
\end{align}
where $\Eeigp = \diag{\mu_{k+1},\ldots,\mu_{2d-k}}$. 
Let $\Em_t=[\nvk_{t}^\top,\B{0}]^\top$,  then the iterates generated by Algorithm~\ref{alg:2} can be expressed as
\begin{align*}
	\Ym_{t+1}
	= 		
	\left(
	\Cm
	\Ym_t
	+
	\Em_t
	\right)
	\Rm_{t+1}^{-1},
\end{align*}
where 
$\Ym_t=[\Wm_t^\top,\Wm_{t-1}^\top]^\top$ with $\Wm_t^\top\Bm\Wm_t=\IM{}$ and $\Wm_{t-1}^\top\Bm\Wm_{t-1}=\Rm_{t}^{-T}\Rm_t^{-1}$. 
Therefore, if we let $\bar{\Ym}_t=\Ym_t(\IM{}+\Rm_{t}^{-T}\Rm_t^{-1})^{-1/2}$, then 
$\bar{\Ym}_t^\top\bigB\, \bar{\Ym}_t=\IM{}.$

By Proposition~\ref{prop:sin:cos:tan}, $\tan\theta\big(\rm{span}(\Ym_{t+1}),\rm{span}(\bar{\Vm})\big)$ w.r.t $\bigB$ is given by
\begin{align}
	\tan\theta\left(\bar{\Ym}_{t+1}, \Vm\right)
	=
		\Vnorm{\Vmp^\top\bigB \,\bar{\Ym}_{t+1}(\Vm^\top\bigB\,\bar{\Ym}_{t+1})^{-1}},
\end{align}
which implies that
\begin{align}
	\tan\theta\left(\bar{\Ym}_{t+1}, \Vm\right)
	&=
		\Vnorm{\Vmp^\top\bigB \,\Ym_{t+1}(\Vm^\top\bigB\,\Ym_{t+1})^{-1}}
	\nonumber\\
	&\hspace{0.0cm}
	=	
	\Vnorm{\Vmp^\top\bigB \,(\Cm\Ym_t+\Em_t)(\Vm^\top\bigB\,(\Cm\Ym_t+\Em_t))^{-1}}
	\nonumber\\
	&\hspace{0.0cm}
	=
		\Vnorm{(\Eeig_\perp\Vmp^\top\bigB\,\Ym_t + \Vmp^\top\bigB\,\Em_t)
		(\Eeig\Vm^\top\bigB\,\Ym_t + \Vm^\top\bigB\,\Em_t)^{-1}}
	\nonumber\\
	&\hspace{0.0cm}
	\leq
		\Vnorm{(\Eeig_\perp\Vmp^\top\bigB\,\Ym_t + \Vmp^\top\bigB\,\Em_t)(\Vm^\top\bigB\,\Ym_t)^{-1}}
		\Big\Vert{\big(\underbrace{\Eeig + \Vm^\top\bigB\,\Em_t(\Vm^\top\bigB\,\Ym_t)^{-1}}_{\B{J}}\big)^{-1}\Big\Vert}
\end{align}
where the third equality follows from~\eqref{eq:C:decomp} and the facts that $\Vm^\top\bigB\,\Vm=\IM{}$ and $\Vmp^\top\bigB\,\Vmp=\IM{}$.
Since
\begin{align*}
	\Vnorm{\B{J}^{-1}}=\frac{1}{\sigma_{\rm{min}}(\B{J})}
	\leq \frac{1}{\sigma_k(\Eeig)-\Vnorm{\Vm^\top\bigB\,\Em_t(\Vm^\top\bigB\,\Ym_t)^{-1}}}
	\leq \frac{1}{\sigma_k(\Eeig)-\Vnorm{\Vm^\top\bigB\,\Em_t}\Vnorm{(\Vm^\top\bigB\,\Ym_t)^{-1}}}
\end{align*}
where $\sigma_k(\cdot)$ denotes the $k$th largest singular value, it follows that
\begin{align}\label{eq:tan:pf}
	\tan\theta\left(\bar{\Ym}_{t+1}, \Vm\right)	
	\leq
		\frac{
			\Vnorm{\Eeig_\perp} \tan\theta\left(\bar{\Ym}_{t}, \Vm\right)
			+
			\Vnorm{\Vmp^\top\bigB\,\Em_t}\Vnorm{(\Vm^\top\bigB\,\Ym_t)^{-1}}
		}{
			\sigma_k(\Eeig)-\Vnorm{\Vm^\top\bigB\,\Em_t}\Vnorm{(\Vm^\top\bigB\,\Ym_t)^{-1}}
		}.
\end{align}
We have 
\begin{align}\label{eq:cos:pf}
	\Vnorm{(\Vm^\top\bigB\,\Ym_t)^{-1}}
	\leq
		\Vnorm{(\Vm^\top\bigB\,\bar{\Ym}_t)^{-1}}\Vnorm{\big(\IM{}+\Rm_{t}^{-T}\Rm_t^{-1}\big)^{-1/2}}
	\leq 
		\Vnorm{(\Vm^\top\bigB\,\bar{\Ym}_t)^{-1}}=1/\cos\theta(\bar{\Ym}_t,\Vm),
\end{align}
where the last step follows from the definition of $\cos\theta(\bar{\Ym}_t,\Vm)$.
Substituting~\eqref{eq:cos:pf} into~\eqref{eq:tan:pf} yields
\begin{align}
	\tan\theta\left(\bar{\Ym}_{t+1}, \Vm\right)	
	&\leq
		\frac{
			\Vnorm{\Eeig_\perp} \tan\theta\left(\bar{\Ym}_{t}, \Vm\right)
			+
			\frac{\Vnorm{\Vmp^\top\bigB\,\Em_t}}{\cos\theta(\bar{\Ym}_t,\Vm)}
		}{
			\sigma_k(\Eeig)-\frac{\Vnorm{\Vm^\top\bigB\,\Em_t}}{\cos\theta(\bar{\Ym}_t,\Vm)}
		}
	\nonumber\\
	&\hspace{0.0cm}
	\leq
		\tan\theta\left(\bar{\Ym}_{t}, \Vm\right)
		\frac{
			\norm{\mu_{k+1}} 
			+
			\frac{\Vnorm{\B{E}_t}_\bigB}{\sin\theta(\bar{\Ym}_t,\Vm)}
		}{
			\norm{\mu_k}-\frac{\Vnorm{\B{E}_t}_\bigB}{\cos\theta(\bar{\Ym}_t,\Vm)}
		},
\end{align} 
where we have used the fact that $\Vnorm{\Vmp^\top\bigB\,\Em_t}\leq\Vnorm{\Em_t}_\bigB=\Vnorm{\B{E}_t}_\Bm$.
As a consequence, if $\B{E}_t$ satisfies 
\begin{align}\label{eq:error:cond}
	\Vnorm{\B{E}_t}_\Bm
	\leq \frac{\norm{\mu_k}-\norm{\mu_{k+1}}}{4}\min\left\{\sin\theta(\bar{\Ym}_t,\Vm),\cos\theta(\bar{\Ym}_t,\Vm)\right\},
\end{align}
then we obtain
\begin{align}\label{eq:tan:ineq}
	\tan\theta\left(\bar{\Ym}_{t+1}, \Vm\right)		
	\leq
		\frac{\norm{\mu_k}+3\norm{\mu_{k+1}}}{3\norm{\mu_{k}}+\norm{\mu_{k+1}}}	\tan\theta\left(\bar{\Ym}_{t}, \Vm\right).
\end{align} 
Having established the inequality~\eqref{eq:tan:ineq}, we can now use similar steps as in the proof of Theorem~\ref{thrm:1} to obtain the following iteration complexity:
\begin{align}\label{eq:iter:complx:topk}
	T=\mathcal{O}\big(\frac{1}{\sqrt{\gap_k}} \log\frac{\tan\angl{}{\bar{\Ym}_0, \Vm}}{\epsilon}\big),
\end{align}
where $\gap_k=1-\eig_{k+1}/\eig_k$.

Similar to the proof of Theorem~\ref{thrm:1}, we define the error in solving the least squares problem as
\begin{align*}
	\textsf{r}\left(\Wm\right) 
	=
		\Mnorm{\Wm-\Bm^{-1}\Am\Wm_t}{\Bm,\rm{F}}^2
	=2\left(f\left(\Wm\right)-f\left(\Bm^{-1}\Am\Wm_t\right)\right).
\end{align*}
Let $\ires=\textsf{r}\left(\Wm_{t}\Zm_t\right)$ and $\dres=\Mnorm{\B{E}_t}{\Bm}^2=\textsf{r}\left(\Wtm_{t+1}\right)$, then it has been shown in~\cite{GJK16} that $\ires\leq 4k\norm{\eig_1}^2\tan\theta(\Wm_t,\Um)$. Combining with~\eqref{eq:error:cond}, it is thus sufficient to decrease the error until
\begin{align}
	\frac{\dres}{\ires}
	\leq 
		\frac{(\norm{\mu_k}-\norm{\mu_{k+1}})^2}{16}
		\frac{
			\min\left\{\sin^2\theta(\bar{\Ym}_t,\Vm),\cos^2\theta(\bar{\Ym}_t,\Vm)\right\}
		}{
			4k\norm{\eig_1}^2\tan^2\theta(\Wm_t,\Um)
		}.
\end{align} 
For $\norm{\eig_{k+1}}\leq 2\sqrt{\beta} < \norm{\eig_{k}}$, we have $\norm{\mu_k} \geq \norm{\eig_k}$, hence if the errors above satisfy 
\begin{align}
	\frac{\dres}{\ires}
	\leq 
		\frac{\gap_k}{128k\gamma^2}
		\min\left\{
			\frac{\sin^2\theta(\bar{\Ym}_t,\Vm)}{\tan^2\theta(\Wm_t,\Um)},
			\frac{\cos^2\theta(\bar{\Ym}_t,\Vm)}{\tan^2\theta(\Wm_t,\Um)}
		\right\},
\end{align} 
where  $\gamma=\norm{\eig_1}/\norm{\eig_k}$, then the condition~\eqref{eq:error:cond} will be satisfied.
Finally, using the same two-phase analysis as in the proof of Theorem~\ref{thrm:1}, we arrive at the desired running time.

\end{document}